\documentclass[letterpaper,12pt]{article}
\usepackage[inner=2.5cm,outer=2.5cm]{geometry}
\usepackage[english,activeacute]{babel}
\usepackage{amssymb,amsmath,alltt}
\usepackage[utf8]{inputenc}
\usepackage{amsfonts}
\usepackage{graphicx}
\usepackage{multirow}
\usepackage{epstopdf}
\usepackage{lscape}
\usepackage{colortbl}
\usepackage{booktabs}
\usepackage{array}
\newcolumntype{L}[1]{>{\raggedright\let\newline\\\arraybackslash\hspace{0pt}}m{#1}}
\newcolumntype{C}[1]{>{\centering\let\newline\\\arraybackslash\hspace{0pt}}m{#1}}
\usepackage{ctable}

\newcommand*{\Resize}[2]{\resizebox{#1}{!}{$#2$}}%

\usepackage{color}
\definecolor{red}{RGB}{250,0,0}
\definecolor{blue}{RGB}{0,0,0}

\newcolumntype{g}{>{\columncolor{Gray}}c}

\definecolor{Gray}{gray}{0.9}
\definecolor{LightCyan}{rgb}{0.88,1,1}

\usepackage[symbol]{footmisc}
\renewcommand{\thefootnote}{\fnsymbol{footnote}}

\usepackage{amsthm}

\newtheorem{algorithm}{Algorithm}

\newcommand{\bq}{\begin{equation}}
\newcommand{\eq}{\end{equation}}

\newcommand{\be}{\begin{eqnarray}}
\newcommand{\ee}{\end{eqnarray}}

\title{A new method to compute periodic orbits in general symplectic maps}
\author{R. Calleja$^\ast$\footnote{calleja@mym.iimas.unam.mx},
        D. del-Castillo-Negrete$^\star$\footnote{delcastillod@ornl.gov},
\setcounter{footnote}{7}
        D. Martínez-del-Río$^\ast$\footnote{david.martinez-del-rio@warwick.ac.uk},
        A. Olvera$^\ast$\footnote{aoc@mym.iimas.unam.mx} \vspace{0.2cm}  \\ 
      \small{$^\ast$ IIMAS-UNAM, Mexico City, Mexico 04510} 
      \phantom{OOOOOOOOOOOOO}  \\   
      \small{$^\star$ Oak Ridge National Laboratory, Oak Ridge, Tennessee 37831-8071}} 
 \date{}

\begin{document}
\maketitle
\renewcommand*{\thefootnote}{\arabic{footnote}}
\setcounter{footnote}{0}

\begin{abstract}
The search of high-order  periodic orbits has been typically restricted to problems with symmetries that help to reduce the dimension of the search space. Well-known examples include reversible maps with symmetry lines. The present work proposes a new method to compute high-order periodic orbits in twist maps without the use of symmetries. The method is a combination of the parameterization method in Fourier space and a Newton-Gauss multiple shooting scheme. The parameterization method has been successfully used in the past to compute quasi-periodic invariant circles. However, this is the first time that this method is used in the context of periodic orbits. Numerical examples are presented showing the accuracy and efficiency of the proposed method.  The method is also applied to verify the renormalization prediction of the residues' convergence at criticality (extensively studied in reversible maps) in the relatively unexplored case of maps without symmetries. 
\end{abstract}


\section{\large{Introduction}}
The systematic study of dynamical systems with low number of degrees of freedom started at the
beginning of last century with the discovery of heteroclinic phenomena in periodic orbits.
Poincaré and Birkhoff studied in detail the existence of periodic orbits in discrete dynamical systems, in particular twist maps that preserve area. The existence and uniformity
of these periodic orbits with respect to their rotation number allowed their use as a tool to analyze other invariant objects that persist in almost integrable maps. From the results obtained by Birkhoff \cite{Meiss92} it was possible to implement a procedure to approximate invariant tori and determine their properties.
Later on, with the arrival of digital computers it was possible to implement search algorithms for these invariant objects through periodic orbits of rotation numbers that approximate the invariant object.
Studies from the 60's and 70's, including Chirikov \cite{Chirikov79}, Greene \cite{Greene79}, Kadanoff \cite{Kadanoff81} and others obtained interesting results relying on numerical experiments performed on increasingly more powerful computers.
In particular, Greene's residue criterion \cite{Greene79}, arguably the most accurate and most used method to determine the persistence or destruction of invariant circles, relies on finding high order periodic orbits for critical parameter values.

As a result from these studies, the renormalization theory of twist maps was introduced
\cite{MacKay82} and it was possible to show the existence of universal behavior for parameters  close to the destruction of invariant curves \cite{Kadanoff81,shenker82}.
All these results were obtained with the help of periodic orbits and, as the numerical capabilities increased, the numerical experiments were performed with periods 
reaching orders of tens of millions and limited only by the arithmetic precision used. 
This opened the possibility of numerically study in detail 
the scenarios predicted by the renormalization theory and the Aubry-Mather theory.

However, the numerical study of periodic orbits in area preserving twist maps has been mostly limited to a  particular kind of maps known as \emph{reversible maps}. The fundamental property of reversible maps is that they can be written as the composition of two involutions which are maps with the property that their composition with themselves is the identity. 
The invariant sets of involution maps are usually one-dimensional sets known as symmetry lines.
Using the theory of DeVogelaere \cite{deVogelaere58,Greene79} it can be proved that two iterates of every \textcolor{blue}{symmetric} periodic orbit lay in these invariant sets. This result allows to simplify the search of periodic orbits with the use of one-dimensional methods (e.g., 1D-quasi-Newton methods). 
The robust and stable  behavior of these 1D methods allows the computation of periodic orbits of very high order, up to $\sim 10^{7}$,  with an accuracy limited only by machine precision 
{\color{blue} Ref.~\cite{Petrov_Olvera08}}. 
This is one of the  reasons why only in reversible maps,
and in particular the  standard map introduced by Chirikov and Taylor \cite{Meiss92, Chirikov79, Taylor69, lichtenberg2013regular}, that it has been possible to study critical dynamical phenomena in detail. 
It is important to point out that reversibility is independent of the twist condition and that there are non-twist maps that are reversible, an extensively studied example being the standard non-twist map Refs.~\cite{del1994dynamics, del_g_m96,del1997renormalization, apte05, fuchss06, gonzalez2014singularity}. Like in the standard map, the use of involutions and symmetry lines 
in the non-twist map have allowed the possibility of studying in great detail the criticality and renormalization of the destruction of shearless invariant circles using periodic orbits of very high order.   

Unfortunately many dynamical systems that can be reduced to area preserving maps can not be studied in the same way as reversible maps either because the lack of knowledge or the actual impossibility of writing them as products of involutions. In this case, the search for periodic orbits must be done in two or more dimensions and, as it is well known, 
{\color{blue}this might compromise the convergence of Newton methods because they typically exhibit 
 very small, in some cases fractal, basins of attraction.}
 This problem is exacerbated in high dimensional maps an example being the 4D Froeschlé map  
 for which the computation of periodic orbits is typically limited to relatively low periods.
 \cite{kook89,olvera94}.

This paper presents a new numerical method to find periodic orbits in general area preserving twist maps without the use of symmetries. 
Among the physics motivations for this study is our recent work on  
self-consistent transport phenomena in a reduced plasma physics model consisting of a large number of standard-like maps coupled by a mean-field \cite{Coreo15, nasm_paper}. The study of the global stability properties of this model required the finding of high order periodic orbits  in 
 nonautonomous maps without known symmetries. 

Our methodology is based on the implementation of the parameterization method developed by de la Llave and Calleja
\cite{Calleja09,Cal-Lla-10b} based on Ref.~\cite{Lla-Gon-Jor-Vil-05} (see also Ref.~\cite{Haro16}), that allows the approximation of continuous invariant objects,  e.g. invariant circles, in two-dimensional twist maps.
 The numerical implementation of the method is based on the computation of Fourier modes of an invariant circle with a fixed rotation number and it provides lower bounds for the critical values of the parameters for which the invariant circle exists as a continuous set. 
 {\color{blue}Among the advantages of this method is that it provides a suitable change of coordinates that allows to reduce the number of operations from $O(N^2)$ to $O(N\log (N))$, where $N$ is the number of points used to represent the invariant curve. This favorable scaling enables the computation of periodic orbits with periods of the order $\sim O(10^7)$.}
The method proposed in the present paper consists of two steps. In the first step, a modified version of the parameterization
method is used to obtain a suitable seed for the use in the second step consisting of a refinement based on the use of a Newton-Gauss method.

The rest of the paper is organized as follows.
Section~\ref{sec:Prelim} presents review material and introduces the rational harmonic map that will be used in the numerical examples. 
Section~\ref{sec:Param_method} discusses the main ingredients of the standard parameterization
method introduced in Refs.~\cite{Lla-Gon-Jor-Vil-05, Cal-Lla-10}.
Section~\ref{sec:Compound} presents the new  compound parameterization method  for the computation of periodic orbits.
Section~\ref{sec_Implementation} discusses the numerical implementation of the proposed method in the context of the standard map and the rational harmonic map including a  comparison with 
the renormalization theory prediction of the residues' convergence at criticality.  
In Section~\ref{sec_Discussion} we present a summary discussion and possible future applications of the new method.

\section{\large{Preliminaries}}
\label{sec:Prelim}

This section reviews basic definitions and concepts, further details can be found in 
standard references including Refs.~\cite{Meiss92,Reichl04,DM19}. 
The maps of interest in the present work are symplectic diffeomorphisms on the cylinder,  $T:\mathbb{S}\times\mathbb{R}\rightarrow\mathbb{S}\times\mathbb{R}$.
A periodic orbit with rotation number $p/q$ on the lift of map $T$, $\tilde{T}:\mathbb{R}^2\rightarrow\mathbb{R}^2$,  satisfies the relation,
\begin{equation}
\tilde{T}^q(z_0) = z_0 + P \,, \,\,\, \qquad \mathrm{where}\,\,\,\,  
 z_0\in\mathbb{R}^2,\,\,
 P= \left( \begin{array}{c} p \\  0 \end{array} \right)\,.
 \label{perio_def}
\end{equation}

A map $T$ is called \emph{reversible} \cite{devaney76,sevryuk86,kook89}
if it can be written as the composition, $T=I_2\circ I_1$, of two functions $I_1$ and $I_2$ with the property,
\begin{equation}
I_k\circ I_k={\rm Id} \,, \qquad k=1,2 \,,
\label{invol_def2}
\end{equation}
where $\rm Id$ is the identity. Functions that fulfill this property are called \emph{involutions}.
A list of properties of reversible maps can be found in Ref.~\cite{roberts94} and further properties of involutions in Refs.~\cite{pina87,roberts92}. 
The full characterization of the general conditions that a map needs to satisfy in order to be written as a product of two involutions is an open problem. In the present paper we say that a map is 
\emph{non-reversible} if there is not a known involution decomposition of the map. 
The invariant sets of involutions are called symmetry lines\footnote{In Ref.~\cite{pina87} a more general definition is used for symmetry lines.}.
{\color{blue}It has been proven that for reversible maps  any symmetric periodic orbit has a point on a symmetry line of the map provided one of the involutions is orientation reversing
\cite{MacKay82,lamb1998time,fox2014critical}.}

The existence of symmetry lines  allows to reduce the search of periodic orbits in
reversible twist \cite{Greene79,Petrov_Olvera08} and non twist   \cite{del_g_m96} maps
from a 2-D problem to a 1-D search on the symmetry lines. 
In 1-D there are plenty of efficient and well behaved methods, including quasi-Newton methods,  to search for the zeros of an equation \cite{dennis96}.
This is why, in reversible maps,  it is possible to find periodic orbits with periods of the order  
$q\sim 10^7$.  As illustrated in Ref.~\cite{Petrov_Olvera08}, the only limit of this type of  computation is the machine precision and not the complexity of the map. 
{\color{blue} However, even when the involutions can be found, the associated
symmetry lines might be quite complex and hard to compute \cite{shi2011reversible}. Even worse, in the case when none of the involutions is orientation reversing there is no guarantee that the symmetry lines will be continuous objects \cite{roberts92} or that they will contain all monotone symmetric periodic orbits.}

\begin{figure}[h!]
     \centering
     \includegraphics[width=0.81\linewidth]{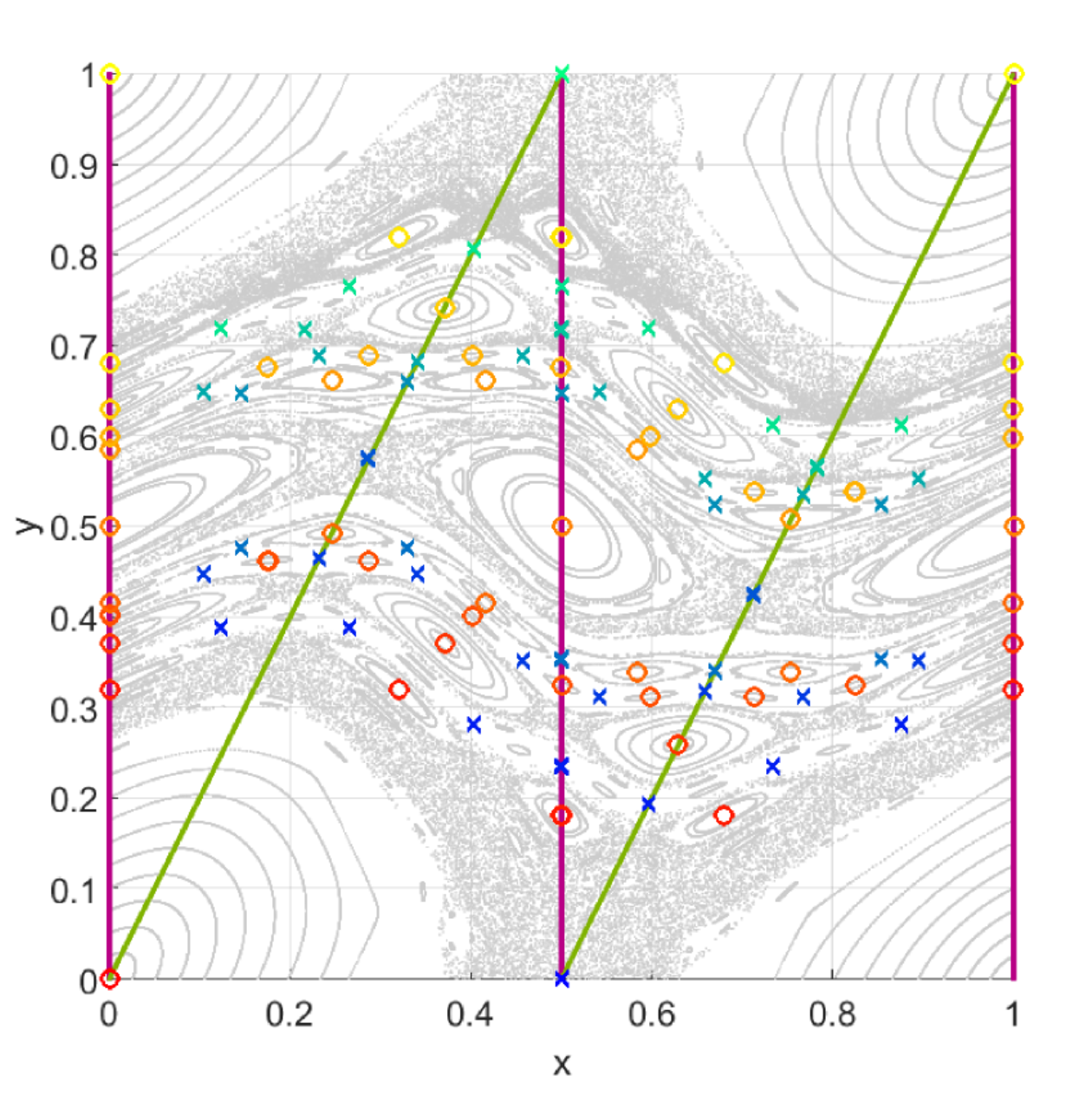} 
     \caption{Example of symmetry lines and periodic orbits in the Chirikov-Taylor map (\ref{stnmap_1}) with $\kappa=0.97$. The purple lines correspond to the symmetry lines of $I_1$ and green lines to symmetry lines of $I_2$. The symmetric periodic orbits search problem is one dimensional since all these orbits have at least one point in one of the symmetry lines.}
     \label{symmlinesfig}
 \end{figure}

As it is well-known, in the case of the \emph{standard map},
\begin{equation}
 \mathcal{S}_\kappa 
 \left( \begin{array}{c}  x\\ y \end{array} \right)  =
    \left(\begin{array}{c}
   {x}+ y -  V'(x) \\ 
 y -  V'(x) 
\end{array} \right)\,,
\label{stnmap_0}
\end{equation}
with %
\begin{equation}
   V'(x) = \frac{\kappa}{2\pi} \sin(2 \pi x) \,.
\label{stnmap_1}
\end{equation}
an involution decomposition, $\mathcal{S}_\kappa=I_2 \circ I_1$, is given by
\begin{equation}
I_1  \left( \begin{array}{c} x \\  y  \end{array} \right)
    = \left( \begin{array}{c} -x \\  y - \frac{\kappa}{2\pi} \sin(2\pi x)  \end{array} \right) \,, \qquad
I_2  \left( \begin{array}{c} x \\  y  \end{array} \right)
    = \left( \begin{array}{c} y - x \\  y   \end{array} \right) \,.
\end{equation}
In this case, as shown in Fig.~\ref{symmlinesfig}, the symmetry lines associated with 
$I_1$ are the vertical lines $x=0$ and $x=1/2$, and those associated with $I_2$ are the lines $x=\frac{y}{2}$ and $x=\frac{y+1}{2}$. Figure~\ref{symmlinesfig} also shows a sample of  elliptic ($\mathsf{O}$) and hyperbolic $(\mathsf{X}$) symmetric periodic orbits which, as indicated above, have two points on each symmetry line.
Note that $x=0$, know as the \emph{dominant} symmetry line \cite{Kadanoff81},
only contains elliptic symmetric periodic orbits. 
%

The linear stability of a periodic orbit is determined by its  \emph{residue}
\begin{equation}
R= \frac{1}{4} [2- \mathrm{Tr}(DT^q)] \,,
\label{greene_res}
\end{equation} 
where $\mathrm{Tr}$ denotes the trace and $D$ the derivative.
Greene's residues criterion \cite{Greene79} establishes a connection between the robustness of a quasiperiodic invariant circle with  rotation number $\omega$ and the behavior of the residues of 
periodic orbits with rotation numbers $p_n/q_n$, such that $p_n/q_n \rightarrow\omega$ as $n\rightarrow \infty$. In particular, if a quasiperiodic invariant circle with a fixed rotation number exists then the residues converge to zero. On the other hand, if the residues diverge, the invariant circle does not exists. The critical case corresponds to the case when the residues converge to a fixed point (or a cycle) different to the $R=0$ trivial integrable limit. 
The work done by Greene on twist  maps was extended to non-twist maps in Refs~\cite{del_g_m96,del1997renormalization,apte05,fuchss06} and other more general cases in Ref.~\cite{Petrov_Olvera08} but always limiting attention to reversible maps with symmetries.

To illustrate and test the proposed method we will use, in addition to the standard map, a map of the form in Eq.~(\ref{stnmap_0}) with perturbation function,  
\begin{equation}
V'(x)= \frac{\kappa}{2\pi} \left(f(x) - \int_0^1 f(s)ds \right) \,,
\qquad \mathrm{where}\,\,\,\,
  f(x)= \frac{\sin(2\pi x+\alpha)}{1-\beta\cos(2\pi x)} \, ,
  \label{anmap_func}
\end{equation}
{\color{blue} with $|\beta|<1$}. We refer to this map as the \emph{rational harmonic} map due to the form of the perturbation function.
Note that for $\beta\neq 0$ the perturbation function $V'(x)$ can be expressed  as an infinite Fourier series and is singular {\color{blue} for $|\beta| \geq 1$}. 
This map has been studied before in Refs.~\cite{Simo_Olvera,Petrov_Olvera08} to test universality results in different kinds of twist maps. 
{\color{blue} 
Following Ref.~\cite{fox2014critical}, it is possible to construct an involution decomposition of this map in the generic case when the perturbation function is neither an odd nor an even function of $x$.} However,
with the exception $\alpha=0$ and $\alpha=\pi$, the associated symmetry lines are not useful for the  computation of periodic orbits and thus this map is a good candidate to test the implementation of the proposed method. 

\section{\large{The parameterization method}}
\label{sec:Param_method}
 The parameterization method was originally introduced by de la Llave et al.\cite{Lla-Gon-Jor-Vil-05} to find an approximate conjugation function between an invariant torus and the rigid rotation over a ideal torus.
The rationale of the method is best understood in the constructive proof of the {KAM} theorem in Ref.~\cite{Lla-Gon-Jor-Vil-05}, which
relies among other things on a Newton iteration in the spirit of Nash-Moser theory, see Ref.~\cite{Zehnder75}. An extra ingredient of the particular implementation described here is related to the area preserving property of the maps of interest  that allows the use of a symplectic change of coordinates, referred in Ref.~\cite{Lla-Gon-Jor-Vil-05} as \emph{automatic reducibility}, that allows to significantly simplify the numerical calculations.

Nash-Moser techniques can be used in algorithms 
that allow to continue smooth functions
$K: \mathbb{S}^1 \to \mathbb{S}^1\times\mathbb{R}$ satisfying the invariance 
equation,
\begin{equation}\label{invariance}
T\circ K(\theta) = K(\theta+\omega) \,,
\end{equation}
where $T:\mathbb{S}^1\times\mathbb{R}\to \mathbb{S}^1\times\mathbb{R}$ is a given twist map and $\omega$ is a Diophantine number, i.e., 
an irrational number such that
\label{diophantine}
for a given $\tau >1$ there is  a constant $\nu > 0$ such that,
\begin{equation}
|\omega\cdot q - p| \geq \nu|q|^{-\tau},\,\,p \in \mathbb{Z},\,\,q\in\mathbb{Z}\setminus{\{0\}} \,.
\end{equation}
 
 Starting from the integrable case of the map $T$, the continuation moves 
a perturbation parameter, e.g. $\kappa$ in the standard map, as close to the breakdown of analyticity of the 
invariant circles as possible. The criterion of breakdown in Ref.~\cite{Cal-Lla-10} is used to determine when the invariant circle ceases to exist. This criterion states that close to the breakdown of analyticity the derivatives of $K$ start to blow up at points of $K(\mathbb{S}^1)$,  in the sense that all the Sobolev $||\cdot||_{H^n}$ norms diverge. Continuation methods like the one presented here have been used in several contexts. See for instance, Ref.~\cite{Calleja09, Cal-Lla-10b},
for models in statistical mechanics, Ref.~\cite{Hug-Lla-Sir-06,
Fig-Luq-Har-15} for examples in symplectic maps, 
Ref.~\cite{Cel-Cal-10, Cal-Fig-12} for conformally symplectic models, 
and Ref.~\cite{Fox-Mei-13} for volume preserving maps.

{\color{blue} The main idea of the parameterization method is to solve the invariance Eq.~(\ref{invariance})  using a perturbative expansion starting from a simplified solution $K_0(\theta)$  with an error\begin{equation}
{e}_0(\theta) = T \circ K_0 (\theta) - K_0(\theta + \omega)\, . 
\end{equation}
The next order of the approximation is written as $K_1(\theta) = K_0(\theta) + \Delta(\theta)$ which by construction has an error 
\begin{equation}
{e}_1(\theta) =DT(K_0(\theta)){\Delta}(\theta) - {\Delta}(\theta + \omega) + {e}_0(\theta)  
+ {\cal O}(\Delta^2(\theta)) \, . 
\end{equation}
From here it follows that if  $\Delta (\theta)$ satisfies
\begin{equation}\label{newton}
 DT(K_0(\theta)){\Delta}(\theta) - {\Delta}(\theta + \omega)
= - {e}_0(\theta) \, , 
\end{equation}
then  second order accuracy, i.e.,  $\|{e}_1\| \sim \|\Delta \|^2 \sim 
\|{e}_0\|^2$, is achieved. 
}

Taking into account the area preserving property of the map $T$, and using the automatic reducibility
property we introduce the change of coordinates ${\Delta}(\theta) =
M_0(\theta) {W}(\theta)$ and, using Eq.~(\ref{newton}), reduce the problem to the solution of two 
cohomological equations
\begin{equation}\label{W2}
W_{(2)}(\theta) - W_{(2)}(\theta + \omega) = -[M_0^{-1}(\theta + \omega) e_0(\theta)]_{(2)}
\end{equation}
\begin {equation}\label{W1}
W_{(1)}(\theta) - W_{(1)}(\theta + \omega) = -[M_0^{-1}(\theta + \omega) {e}_0(\theta)]_{(1)} - 
S_0(\theta) W_{(2)}(\theta) \, ,
\end{equation}
where the subindices $_{(1)}$ and $_{(2)}$ denote the first and second components of the vector, 
\begin{equation}\label{Mmatrix}
M_0(\theta) 
  = \left( \begin{array}{c} DK_0(\theta)   \\ \end{array} 
     \Big| \begin{array}{c} J^{-1}DK_0(\theta){N}_0(\theta) \\  \end{array}  \right) \,, 
     \qquad {N}_0(\theta) = [DK_0(\theta)^t DK_0(\theta)]^{-1} \, , 
\end{equation}
where $J$ is the symplectic matrix, and
\begin{equation}
S_0(\theta) = N_0(\theta + \omega) DK_0^t(\theta + \omega) 
DT(K_0(\theta)) DK_0(\theta) N_0(\theta) \,.
\end{equation}
It is worth mentioning that the function $S_0(\theta)$ is related to the local twist condition on the invariant circle $K_0$.

A necessary condition to solve the functional system of equations (\ref{W2})-(\ref{W1}) in Fourier space   is that the average of the right hand side of both equations vanishes, or is of order  ${O}(\|{e}_0\|^2)$. This can be done easily for (\ref{W1}) since, if the map is twist, the average of $W_2$ is a free parameter that can be adjusted as needed\footnote{If the map were not twist, it could be the case that $S_0(\theta)=0$.}. However for the right hand side of (\ref{W2}) from  \emph{Lemma 9} from Ref.~\cite{Lla-Gon-Jor-Vil-05},
\begin{equation}
\int_{\mathbb{S}^1}[M_0^{-1}(\theta + \omega) {e}_0(\theta)]_{(2)}\,d\theta = {O}(\|{e}_0\|^2)\,.
\label{lemma9_Lla05}
\end{equation}
\textcolor{blue}{Note that if $K(\theta)$ solves Eq.~(\ref{invariance}), then 
$K(\theta+\sigma)$ also satisfies Eq.~(\ref{invariance}) for any $\sigma$.}
Further details on the particular implementation of the method used can be found in  Ref.~\cite{nasm_paper}.

\section{\large{A new compound method}}
\label{sec:Compound}

The proposed compound method consists of the composition of two methods: a modified parameterization method, described in Secs.~\ref{subsec:modified} and 
\ref{subsec:tracking}, and a 
Newton-Gauss method,  described in Sec.~\ref{subsec:Newton2D}.  
\subsection{Modified parameterization method}
\label{subsec:modified}
{\color{blue}
The main advantage of the parameterization method, and the reason why we use it,  is that it does not rely on the use of symmetries.  However, the original version of this method was used to study quasiperiodic invariant circles and one of our contributions is to extend its application to describe monotone periodic orbits in twist maps, which by construction are discrete (non-continuous) objects guaranteed to exist 
by Birkhoff's theorem and the Aubry-Mather theory.
As explained in detail below, our idea is to use the parameterization method to approximate a curve that contains the periodic orbit and satisfies an invariance equation at these points. 
By doing this, we are effectively reducing the search of the periodic orbit to a lower dimensional object.}

 When one naively attempts to apply the parameterization method to a Liouville rotation number the computation will not converge due to the ``small denominators" problem, or the existence of 
 ``zero denominators" in the case of rational rotation numbers. 
 These zero denominators appear \textcolor{blue}{when trying to solve Eqs.~\eqref{W2} 
and \eqref{W1} that have the generic form
\begin{equation}
\varphi(\theta) - \varphi(\theta + \omega) = \eta(\theta) \,,
\label{reduced_inv}
\end{equation}
where $\eta$ has zero average.
 The zero denominators are due to the fact that we are looking for a
   function $\varphi$ that
does not satisfy equation \eqref{reduced_inv} for every $\theta \in \mathbb{S}^1$, but only on a 
uniform grid of points $\{\theta_j\}_{j=1}^{q} \in \mathbb{S}^1$.
This will be enough since the invariance condition (\ref{invariance}) only needs to be
valid at points on the periodic orbits.
It is worth noticing that the cohomological operator on the l.h.s. of \eqref{reduced_inv} has
  a null space when $\omega = p/q$ so the solution is not unique. 
 This is consistent with the fact that there are infinitely many curves
that satisfy the functional equation
(\ref{reduced_inv}) at the grid of points corresponding to the periodic orbit.
We will approximate a solution to the invariance equation
 by choosing one of the approximate solutions $\varphi$ of the cohomology equation \eqref{reduced_inv}.}

\textcolor{blue}{In Fourier space, the solution of equation (\ref{reduced_inv}) is given by
\begin{equation}\label{cohom_sol}
  a_k = \frac{b_k}{1-e^{2\pi i\omega k}}\,, \qquad \mathrm{if}\,\, e^{2\pi i\omega k}\neq 1\, 
\end{equation}
with
 \begin{equation}
  \varphi(\theta) = \sum_{k=-\infty}^{\infty} a_k e^{2\pi ik\theta} \,, \qquad
  \eta(\theta) = \sum_{k=-\infty}^{\infty} b_k e^{2\pi ik\theta} \,,
 \end{equation}
 If $\omega = p/q$, then obviously $ e^{2\pi i\omega k}= 1$ for $k = m q$ with $m \in \mathbb{Z}$. To approximate the solution to equations \eqref{W1} and \eqref{W2}, we use the fact that the Fourier coefficients 
$a_k$ can be determined as long as $k\neq mq$. We approximate the solution to equation \eqref{reduced_inv} in two steps. 
 First, we look for a solution \eqref{cohom_sol} in the subspace, $a_{mq}=0$ for $m\in\mathbb{Z}$.
Then we adjust the phase $\sigma$ to minimize the error 
\begin{equation}
{e}_1(\theta_j + \sigma) = T \circ K_1 (\theta_j + \sigma) - K_1(\theta_j + \sigma + \omega)\,
 \end{equation}
 at the grid of points corresponding to the periodic orbit.}

Using this solution as the new $K_0$ we repeat the perturbative computation of the invariance described above up until we reach a prescribed error tolerance, $\tilde \varepsilon$.
In practice, we found that even though
convergence of the method is not yet proven, it provides a good seed for the 
rapid convergence of Newton-Gauss method that will be explained in Section \ref{subsec:Newton2D}.

The number of Fourier modes used in the numerical computations is always
finite. However it is guaranteed by the theory behind the parameterization
method that for a fixed upper bound for the error of the approximate solution,
$\|e_1\|<\tilde{\varepsilon}$ there is an optimal maximum number of
Fourier coefficients $N^*(\tilde{\varepsilon})$ to be considered so the norm of
the remaining tail of harmonics can always be safely included inside the error
of the approximate solution \cite{Cal-Lla-10}. 
\textcolor{blue}{We conjecture that}
the convergence of the new method 
\textcolor{blue}{can be proven, although at the present time}
there is no clear argument of how many harmonics 
are necessary 
to use. If the number of harmonics $N^*$ is fixed a priori, then the method, when it converges, should give an optimal trigonometric polynomial of degree $N^*$ such that the parametric curve pass close to the points of the periodic orbit. And even if the method does not converge and gives an approximate value of a point of the periodic orbit with a bounded error, it may suffice to be used as seed for {\color{blue} the Newton-Gauss method in Section~\ref{subsec:Newton2D}}. So a convergence theorem is not strictly needed for the successful implementation of the method.
We compared the efficiency of the method for different number $N^*$ of harmonics in the approximation. For small values of $q$, the total number of harmonics $N^*$ was taken as $2q$, while for larger values $N^*\sim 4q$.  Increasing the value of
$N^*$ above $4q$ did not reflect in major improvements to the method. 

Although the hypothesis of the existence of a parametric curve $K$ that crosses
all the points of the periodic orbits of rotation number $p/q$,
which will be denoted as $K_{p/q}$, does not seem
unreasonable for maps like the standard map (\ref{stnmap_1})
and probably any locally twist map, it is an open problem to prove its existence. The numerical evidence indicates this to be the case and it should be possible to find a mathematical proof, in the line of Ref.~\cite{JGlz20}. %
However, the present work  focuses on the implementation of the hybrid method and on presenting numerical evidence.
If an existence theorem is proven or if its error can be bounded, an interesting question would be to relate properties such as the regularity of the parametric curves $K_{p_n/q_n}$ to a 
$K_\omega$, when $\left\lbrace p_n/q_n \right\rbrace_{n\in\mathbb{N}}$ is a sequence that converges to a Diophantine rotation number $\omega$. A different approach could be to use an a posteriori theorem of the type in Refs.~\cite{hungria2016rigorous, castelli2017parameterization, burgos2019spatial} 
to prove the existence and uniqueness of these periodic points.
 
\subsection{Phase tracking}
\label{subsec:tracking}
The modified parameterization method is expected to yield a parametric curve $K_{p/q}$ that approximately contains the periodic orbit of interest. The problem now is to find the correct phase $\sigma$, for a point $\theta_j$ in the grid since all the others are related by $\theta_{j+1}=\theta_j +p/q$. 
For simplicity we chose the $\theta_j$ closest to zero, define 
$\tilde{\theta} \equiv \theta_j+\sigma \in [0,1/q)$, and write the periodic orbit condition in
Eq.~(\ref{perio_def}) in terms of $K_{p/q}$ as
\begin{equation}
T^q(K_{p/q}(\tilde{\theta})) = K_{p/q}(\tilde{\theta}) +P  \, ,
\label{FK_phase}
\end{equation}
where $P=(p,0)^t$,
and the invariance condition in Eq.~(\ref{invariance}) for the periodic orbit as 
\begin{equation}
T\circ K_{p/q}(\tilde{\theta}) = K_{p/q}\circ \mathsf{R}_{p/q}(\tilde{\theta}) \,,
\label{conjugation_KRT}
\end{equation}
where $\mathsf{R}_{p/q} (\tilde \theta)=\tilde{\theta}+ p/q$. 
Note that any shift of the form $\tilde \theta \rightarrow \tilde \theta+n p/q$, with $n\in \mathbb Z$,  will also be a solution of Eqs.~(\ref{FK_phase}) and 
(\ref{conjugation_KRT}).

\begin{figure}[h!]
     \centering
     \includegraphics[width=0.95\linewidth]{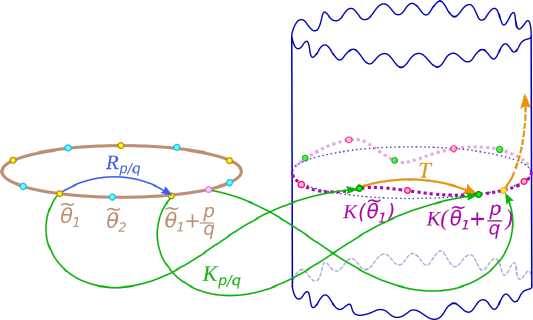} 
     \caption{Illustration of the components of the modified parameterization method, for which the parameterization $K_{p/q}$ is only dynamically consistent with the map $T$ for periodic orbits points.}
     \label{param_fig_po}
 \end{figure}

Usually, to find the
complete set of  elliptic and hyperbolic periodic orbits we only need to find two
independent $\tilde{\theta}_1$ and $\tilde{\theta}_2$ where by independent we mean that there is no $n\in\mathbb{Z}$
such that $\tilde{\theta}_2=\tilde{\theta}_1+np/q\,\,\mathrm{mod}\,2\pi$,
see Fig.~\ref{param_fig_po}.
{\color{blue}The
parameterization method applied to periodic orbits can be considered as an
implementation of the flux-minimizing curves calculation for periodic orbits
described in Refs.~\cite{dewar1992flux} and \cite{dewar2012action}.}
%

There are two continuous functions that can be used to  guide the correction of the  phases $\tilde{\theta}_{1}$ and $\tilde{\theta}_{2}$ for elliptic and hyperbolic periodic points:
the error of the periodic orbit and the residue.
The error can be evaluated at one point (fixed $\theta$)  of the orbit,
\begin{equation}
E(\theta)=||T^q(K_{p/q}(\theta))-K_{p/q}(\theta)-P|| \,, 
\label{true_error}
\end{equation}
or computed as an average over the complete cycle of the periodic orbit,
\begin{equation}
 \tilde{E}(\theta) = 
               \frac{1}{q} \sum_{n=0}^{q-1} 
	    \left| T \circ K_{p/q}(\theta+n p/q) - K_{p/q}[\theta + (n+1) p/q] \right|  \,. \\ \label{error_sob_1}
\end{equation}
Although both functions are regular, $ \tilde{E}$ is better behaved for small values and thus we will use it in the calculations. 
On the other hand, the residue can be defined as a continuous function of $\theta$ using
\begin{equation}
R(\theta)= \frac{1}{4} \left\{2- \mathrm{Tr}\left[DT^q\left (K_{p/q}\left(\theta\right)\right)\right]\right\}   \label{true_residue}
\, . 
\end{equation}
The functions in (\ref{error_sob_1}) and (\ref{true_residue}) are defined on $\theta\in[0,2\pi)$. However, to find $\tilde{\theta}_i$ it is enough to restrict the search to the interval $\theta\in[0,1/q)$. 
Our ansatz  is that the minima of the error (either $E(\theta)$ or $\tilde{E}(\theta)$) correspond to the points in the
curve $K_{p/q}(\theta)$ that are closer to the actual periodic points,
and  the extrema of $R(\theta)$ indicate the location of the elliptic and hyperbolic periodic orbits. 
{\color{blue}However,  care must be taken in the computation of %
$R(\theta)$ 
 because the hyperbolic behavior of the invariant manifolds can cause the errors to grow exponentially in the neighborhood of periodic orbits \cite{Llave_Olvera06}}.
  
The underlying idea behind the parameterization method is that
 the exact dynamic of the map $T$ is conjugated to a rigid rotation over a
 parametric curve $K$. So it is reasonable to expect 
 that for $\theta\in(\tilde{\theta}_i -\delta,\tilde{\theta}_i +\delta)$ 
 the dynamics of $T(K_{p/q}(\theta))$ can be
 approximated by the conjugation, 
\begin{equation}
T^n(K_{p/q}(\theta)) \approx K_{p/q}(\theta+np/q)  \, .
\label{aprox_invariance}
\end{equation}
Then, using 
Eq.~(\ref{aprox_invariance}) in Eq.~(\ref{true_residue})  gives the following expression for the residue         
\begin{equation}
\tilde{R}(\theta)= \frac{1}{4} \Big[ 2- \mathrm{Tr}\Big(\tilde{\mathcal{M}}^q(\theta)\Big) \Big] \,,
\label{aprox_res}
\end{equation}
where 
\begin{equation}
\tilde{\mathcal{M}}^q(\theta) = DT[K_{p/q}(\theta+(q-1)p/q)]
                       \cdots DT[K_{p/q}(\theta)]  \, .
                       \label{DT_tilde}
\end{equation}
The advantage of using Eq.~(\ref{aprox_res})  and Eq.~(\ref{error_sob_1})
for the numerical computations  is that $\tilde E(\theta)$ and $\tilde R(\theta)$ are 
well-behaved because they are  evaluated for  smooth functions over bounded domains. 

The standard numerical computation of the residue using the iteration of the map and its Jacobian  can diverge for small errors in the approximation of the periodic orbit, giving values of the order $\sim 10^{50}$ for periods of the order  $\sim 10^{2}$. In contrast, the  computation using the conjugation function $K_{p/q}$ keeps, by construction, the value of the residue always bounded. 
For example, for periods of the order $\sim10^2$ and an error in the periodic orbit of the order
$\sim 10^{-5}$ it is possible to accurately compute residues with values of the order  $\sim 10^{-15}$.
{\color{blue} 
We conjecture that a computer assisted proof of the convergence of the modified parameterization method could be based on  the use of Newton-Kantorovich theorems in
radii polynomials as done in Ref.~\cite{burgos2019spatial}. However, further exploring this interesting issue is outside the scope of this work.} 
%
 \begin{figure}[h!]
     \centering
     \includegraphics[width=0.92\linewidth]{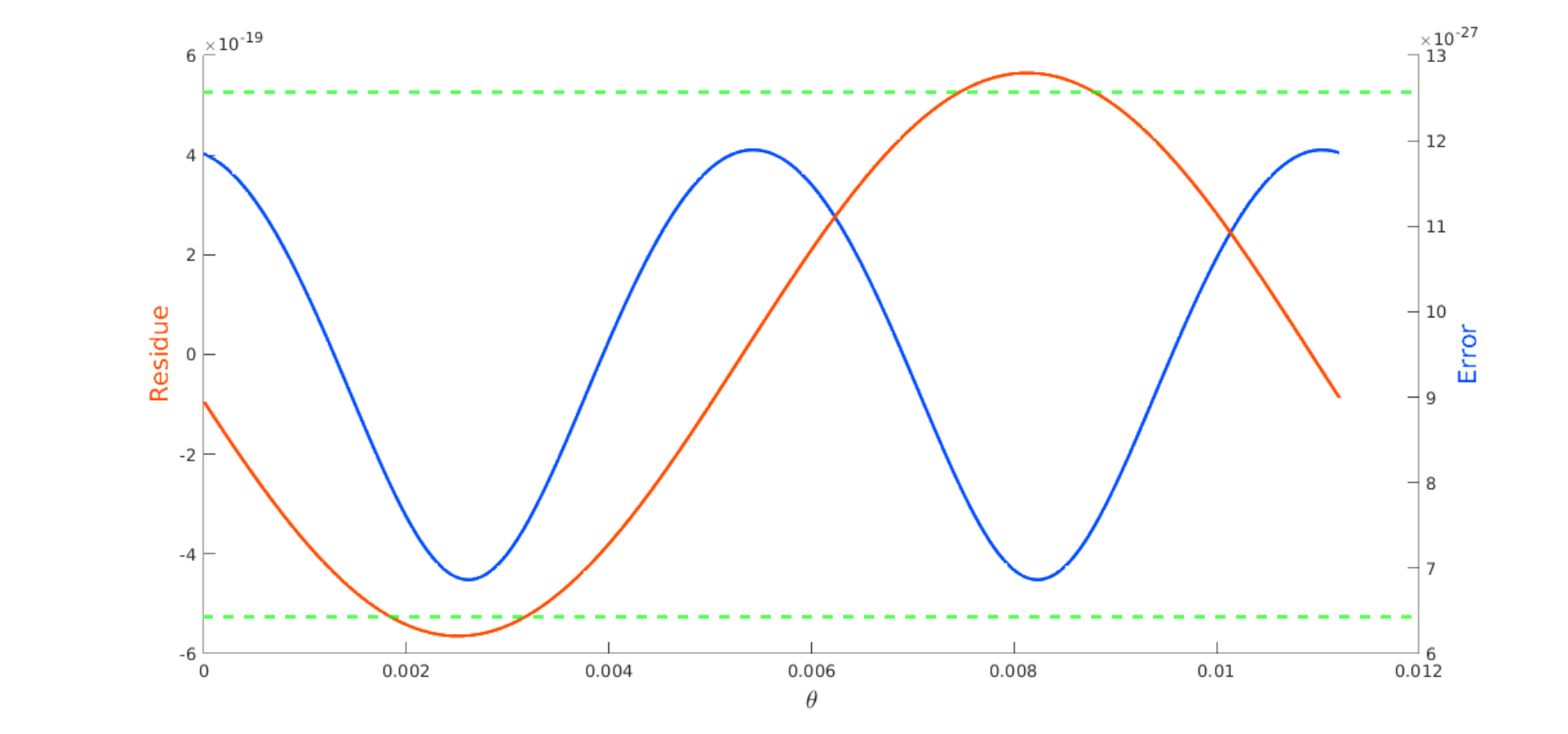}
     \includegraphics[width=0.92\linewidth]{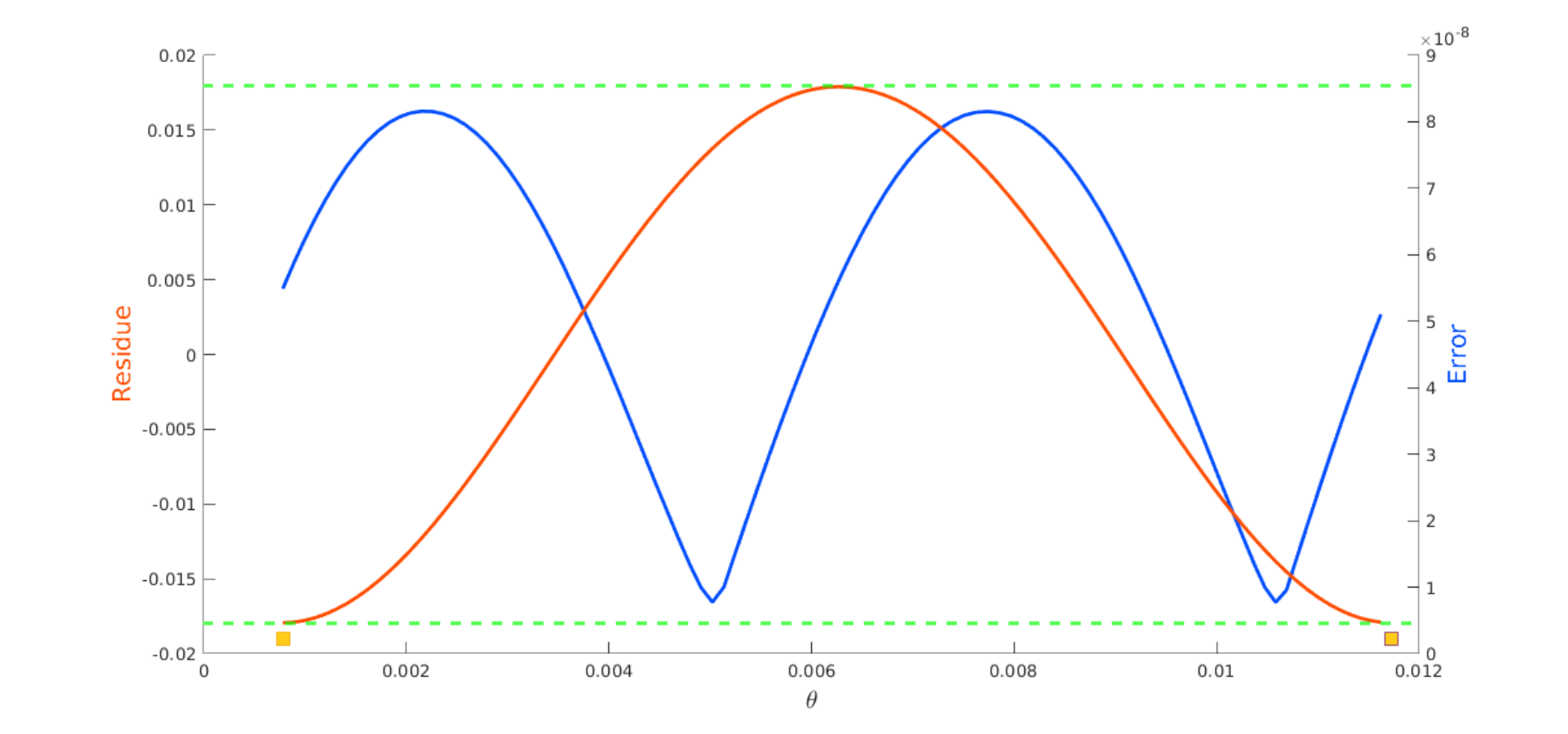} 
     \caption{\textcolor{blue}{Error $\tilde{E}$, Eq.~(\ref{error_sob_1}), (blue) and  
      residue $\tilde{R}$, Eq.~(\ref{aprox_res}), (red) of the approximated        $55/89$-periodic orbit of the \emph{rational harmonic map} (\ref{anmap_func})
     with $(\alpha,\beta)=(2.5,0.37)$ and $\kappa=0.4$ (top panel), $\kappa=0.9$ (bottom panel).     
     The green dashed horizontal lines correspond to the values of the residues  computed using the standard direct Newton method:
   top panel $-5.253\times 10^{-19}$ (hyperbolic) and $5.254\times 10^{-19}$ (elliptic), 
   bottom panel $-0.01797122$ (hyperbolic) and $0.01797122$ (elliptic).  
   The orange squares in the bottom panel indicate the angular position of the hyperbolic periodic points computed with the compound method discussed in Sec.~\ref{subsec:compound}.}}
     \label{error_res_fig}
 \end{figure}

The regularity of the function $\tilde{R}$ and the fact that it should be positive (negative) close to the elliptic (hyperbolic) periodic points implies that, as shown in Fig.~\ref{error_res_fig}, $\tilde{R}$ is an oscillatory function with mean close to zero. \textcolor{blue}{In the top panel of Fig.~\ref{error_res_fig} the reader can notice that the green dashed lines corresponding to the reference residue values (computed using the standard method) of the elliptic and hyperbolic periodic points do not coincide with the minimum and maximum of $\tilde{R}(\theta)$. This can be explained by the error in the minima of the parameterized curve $K_{p/q}(\theta)$ ($>6\times10^{-25}$) that informs us that the best points in the curve to approximate the periodic points are not yet close enough to be considered a good approximations.} 
\textcolor{blue}{In contrast, in the bottom panel of Fig.~\ref{error_res_fig}, which corresponds to parameter values closer to critical, we can observe that there is an appreciable shift between the error minima and the residue critical points, while the last ones are much closer to the correct angular position of the periodic orbit points.}

In all the cases considered it was observed that, like in Fig.~\ref{error_res_fig}, 
 the minima $\{\theta_n^{({min})}\}_{n=1,...,q}$ of $\tilde{E}$  %
are always in the vicinity of the maxima and minima $\{\theta_n^{(crit)}\}_{n=1,...,q}$ of $\tilde{R}$. 
Although there is not yet a  mathematical proof of this observation,  the extrema of $\tilde{R}$  seem to be a useful estimator 
of the position $\tilde{\theta}$ of the periodic points in
$K_{p/q}(\tilde{\theta})$. 
{\color{blue} Figure~\ref{error_res_fig} also shows the dependence of this result on
$\kappa$. In particular, the coincidence of the minima of $\tilde{E}(\theta)$ and the extrema of $\tilde{R}(\theta)$ deteriorates as $\kappa$ increases from $0.4$ to $0.9$. In particular, in the bottom panel of Fig. \ref{error_res_fig} corresponding to $\kappa=0.9$ it can be observed that the extrema of $\tilde{R}(\theta)$ continue to be in the neighborhood of the periodic points while the minima of $\tilde{E}(\theta)$ are shifted away. 
Most likely this mismatch is due to the fact that, close to criticality, the number of Fourier modes 
significantly increase due to the eventual loss of regularity. 
Also note that the computation of the residue is based on the linearized dynamics whereas the computation of the error takes into account the full dynamics. 
Because of these observations, we opted to use as an indicator of the angular position of the periodic orbit points the extrema from $\tilde{R}(\theta)$ and then use them as seed data for the Newton-Gauss method described below.}

\subsection{Newton-Gauss method}
\label{subsec:Newton2D}

In this section we describe the Newton-Gauss method which is a high precision Newton-like method developed by \`A. Haro and collaborators \cite{figueras2015different}
(see also Ref.~\cite{kook89app}) to calculate (or refine) periodic orbits of high order.

Let  $\tilde{\mathcal{Z}}$ be an approximation of the periodic  orbit $\mathcal{Z}$. Then, on the lift, the points $\tilde{z}_{i}\in\tilde{\mathcal{Z}}$ satisfy 
\begin{equation}
\left\lbrace \begin{array}{ll} 
   T(\tilde{z}_{q-1}) - \tilde{z}_0 - P \hspace{0.38cm} = \mathsf{e}_0\,, \\
   T(\tilde{z}_0) - \tilde{z}_1  \hspace{1.63cm} = \mathsf{e}_1 \,, \\
   \hspace{2.85cm}  \vdots & \\ 
   T(\tilde{z}_{q-2}) - \tilde{z}_{q-1}  \hspace{0.9cm} = \mathsf{e}_{q-1} \, ,
       \end{array} \right.
       \label{multi_shooting}
\end{equation}
where $\mathsf{e}_i\in\mathbb{R}^2$, $\|\mathsf{e}_i\|   \ll 1$ for $i=0,1,\dots,q-1$ and $P$ is defined in (\ref{perio_def}).
The error in the approximation is given by 
\bq
E=\|T^q(\tilde{z}_0)-\tilde{z}_0 -P\|\, .
\label{E_NG}
\eq

To compute the periodic orbit $\mathcal{Z}$ we need to find the roots of (\ref{multi_shooting}), 
 i.e. $\left\lbrace \tilde{z}_0,\tilde{z}_1, \dots,\tilde{z}_{q-1}\right\rbrace$ such that $\mathsf{e}_k=0$. 
 Denoting 
 (\ref{multi_shooting}) as the 
 $\mathbb{R}^{2q}\rightarrow\mathbb{R}^{2q}$ function 
 $G(\tilde{\mathbf{z}})=\mathbf{e}$
where $\tilde{\mathbf{z}}=(\tilde z_0,\dots, \tilde z_{q-1})^t$
and 
$\mathbf{e}=(\mathsf{e}_1,\dots,\mathsf{e}_{q-1})^t$, 
the Newton step is given by the contracting map,
\begin{equation}
\hat{\mathbf{z}} = \tilde{\mathbf{z}} - DG^{-1}(\tilde{\mathbf{z}})\cdot 
G(\tilde{\mathbf{z}})\,,
\label{multi_newton}
\end{equation}
where $\tilde{\mathbf{z}}$ is assumed to be inside an open ball sufficient small of a root $\mathbf{z}_*$ of $G(\mathbf{z}_*)=0$.
Defining $\Delta=\hat{\mathbf{z}}-\tilde{\mathbf{z}}$ 
Eq.~(\ref{multi_newton}) can be written as,
\begin{equation}
DG(\tilde{\mathbf{z}}) \Delta = -\mathbf{e}\,,
\label{multi_newton2}
\end{equation}
where $DG(\tilde{\mathbf{z}})$ is a $2q\times2q$ matrix. In more detail (\ref{multi_newton2}) has the form,

\begin{equation}
\left( \begin{array}{cccccc} 
   -I       &\bigcirc &\bigcirc  &\ldots &\bigcirc &DT(z_{q-1})\\
   DT(z_0)  & -I      &\bigcirc  &\ldots &\bigcirc &\bigcirc\\
   \bigcirc &DT(z_1)  &-I        &\ldots &\bigcirc &\bigcirc\\
   \vdots   &\vdots   &\ldots    &\ddots &\vdots   &\vdots\\
   \bigcirc &\bigcirc &\bigcirc  &\ldots &DT(z_{q-1})  &-I
       \end{array} \right)
\left( \begin{array}{c} 
   \Delta z_0 \\ \Delta z_1 \\ \Delta z_2 \\ \vdots \\ \Delta z_{q-1} 
       \end{array} \right)
= \left( \begin{array}{c} 
       -\mathsf{e}_0 \\ -\mathsf{e}_1 \\ -\mathsf{e}_2 \\ \vdots \\ 
       -\mathsf{e}_{q-1} 
       \end{array} \right) \,,
\label{multi_Jacobean}
\end{equation} 
where $I$ is the $2\times2$ identity matrix, $\bigcirc$ the $2\times 2$ zero matrix and $DT(z_j)$ is       the Jacobian matrix of the map $T$  evaluated at $z_j$.
{\color{blue} 
Note that $DG(\tilde{\bf z})$ is sparse since it only has non zero entries in blocks along two diagonals and in the last column.}

\begin{figure}[h!]
     \centering
     \includegraphics[scale=1.0,trim={1.5cm 1.1cm 1.1cm 1.0cm}]{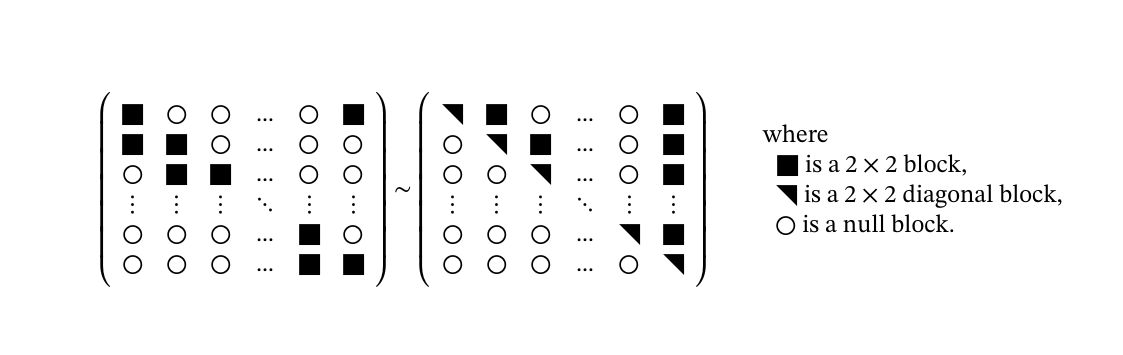}
     \caption{Block representation of matrix in Eq.~(\ref{multi_Jacobean}) before and after triangularization.}
     \label{newton_gauss_r0}
 \end{figure}

To solve the linear system in Eq.~(\ref{multi_Jacobean})
we use a Gauss elimination method on the columns which, 
as shown in Fig.~\ref{newton_gauss_r0},
 transforms  $DG(\tilde{\mathbf{z}})$ on an upper triangular matrix.
 %
The memory storage of this procedure is proportional to the orbit of period $q$. In particular, we need $12q+4$ memory locations to store the $DG(\tilde{\mathbf{z}})$ matrix and $4q$ locations to store the vectors $\Delta$ and $\mathbf{e}$.

Note that the Newton-Gauss method just described is different from the 
standard Newton method based on the iteration
\begin{equation}
 z_{i+1} = z_i - \left[(DT^q)^{-1}(z_i) \right] \mathcal{G}(z_i)\,
 \label{newton_2D}
\end{equation}
where  $\mathcal{G}(z) = T^q(z) - z - P=0$.
Although  Eq.(\ref{newton_2D}) seems simpler than Eq.(\ref{multi_newton2}), since the dimension  of  $z$ is just two (rather than $2q$), Eq.~(\ref{newton_2D}) is known to 
be highly unstable in two or higher dimensions, and it is not reliable for finding periodic orbits of period higher than $10^3$.
On the other hand, the Newton-Gauss method is very stable and allows to find periodic orbits of period $\lesssim 10^7$ with the use of quadruple arithmetic precision. 
The robustness of the method relies on the fact that it solves (by Gauss method) one iteration of the Newton step for one iteration of the map on every point of the orbit at the same time. This process is very fast and computationally cheap in terms of the required memory storage.
Also, it is important to point out that the Newton-Gauss method can be used as a continuation method by considering small variation of the parameters.

\subsection{The compound method}
\label{subsec:compound}
As summarized in the algorithm below, 
the proposed compound method consists of two steps. The first step is based on the use of the modified parameterization method in Secs.~\ref{subsec:modified} and \ref{subsec:tracking} to compute  an approximate periodic orbit. In the second step, this approximate solution is used as the seed of the Newton-Gauss method in 
Sec.~\ref{subsec:Newton2D} to accurately compute the periodic orbit. 

%
\begin{algorithm}
  \label{algoritmo3}
  \begin{itemize}
    \item[]
      \item[I] Compute seed using the modified parameterization method 
  \begin{itemize}
  \item[]
  \item[1)] Given $K_0$, let $e_0 (\theta) = {T} \circ K_0(\theta) - K_0(\theta + p/q)$
  \item[2)] Compute the matrix $M_0(\theta)$ in Eq.~\eqref{Mmatrix}
  \item[3)] Get  $W_2(\theta)$ by solving Eq.~\eqref{W2} eliminating all resonant terms,
  $c_{nq} e^{\pm i nq \tilde{\theta}}$, in the Fourier series of $e_0$ and $W_2$
    \item[4)] Choose the average, $\int_{\mathbb{S}^1} W_2(\theta)d\theta$,
   to guarantee that the average of $[M_0^{-1}(\theta + \omega) e_0(\theta)]_1 + S_0(\theta) W_2(\theta)$ 
   is at least $O(\|e_0\|^2)$ 
     \item[5)] Solve for $W_1(\theta)$ from Eq.\eqref{W1}, eliminating all the resonant terms.
  \item[6)] Compute the step $\Delta(\theta) = M_0(\theta) W(\theta)$, and construct the next order of the approximation $K_1(\theta) = K_0(\theta) + \Delta(\theta)$
  \item[8)] Set $K_0(\theta) = K_1(\theta)$ and go to step 1) until an a priori fixed bound,
  ${\rm min} \{ |e_0 (\theta)| \} < \tilde{\varepsilon}$, is satisfied
  \item[9)] Find two adjacent local minimum and maximum $\left\lbrace \tilde{\theta}_1,\tilde{\theta}_2 \right\rbrace$ of $\tilde R(\theta)$
    \end{itemize}
          \item[II]  Refine seed using the Newton-Gauss method
   \begin{itemize}
  \item[10)] Apply the Newton-Gauss method Eqs.~(\ref{multi_newton})-(\ref{multi_Jacobean}) to $(x_0,y_0)=K_0(\tilde{\theta}_k)$, $k=1,2$ to obtain a point from each one of the hyperbolic and elliptic periodic orbits with the required precision.
  \end{itemize}
    \end{itemize}
\end{algorithm}
%

\section{\large{Numerical implementation of the compound method}}
\label{sec_Implementation}

In this section we discuss the application of the compound method for the computation of periodic orbits in the 
Chirikov-Taylor standard map in Eq.~(\ref{stnmap_1})  and the rational harmonic map in Eq.~(\ref{anmap_func}).
Although the periodic orbits in the standard map can be computed very efficiently using symmetry lines, this map is a good starting point to illustrate the implementation and test the accuracy of the new compound method against well-established numerical results.  
Conceptually speaking, the rational harmonic map presents a more interesting application because it is a much less studied map with three free parameters lacking useful symmetries
that could reduce the computation of periodic orbits to a one-dimensional search.

\subsection{Standard map}

The first test of the \emph{compound method} was done using the Chirikov-Taylor
map Eq.~(\ref{stnmap_1}) to compare with the periodic orbits obtained via symmetry lines
algorithms.  This test allowed to check the implementation of the algorithm, 
obtain bounds for the free parameters of
the method (including the number of Fourier modes and the continuation step size and others)
and compare the execution run times and the 
numerical precision achieved. 

To test the capabilities of the compound method as a continuation method to obtain periodic orbits, first we tested the modified parameterization method alone (without the Newton-Gauss refinement) on the Chirikov-Taylor map.  
The continuation started from the integrable case, $\kappa=0$, to the neighborhood of the
critical parameter value found by Greene, $\kappa_c=0.971635046$, for which all the periodic orbits with rotation number that approximate the inverse of the golden mean have the same critical residue value, $R=-0.255426$ (for hyperbolic orbits) and $R=0.2500888$ (for elliptic orbits) \cite{Reichl04}.
From earlier tests of the code it was found the need to perform the computations with quadruple precision ($O\sim 10^{-30}$), for both the computation of the seed using the modified parameterization (in particular in the computation of the FFT) and the Newton-Gauss refinement.
%
 \begin{figure}[h!]
     \centering
     \includegraphics[width=0.98\linewidth]{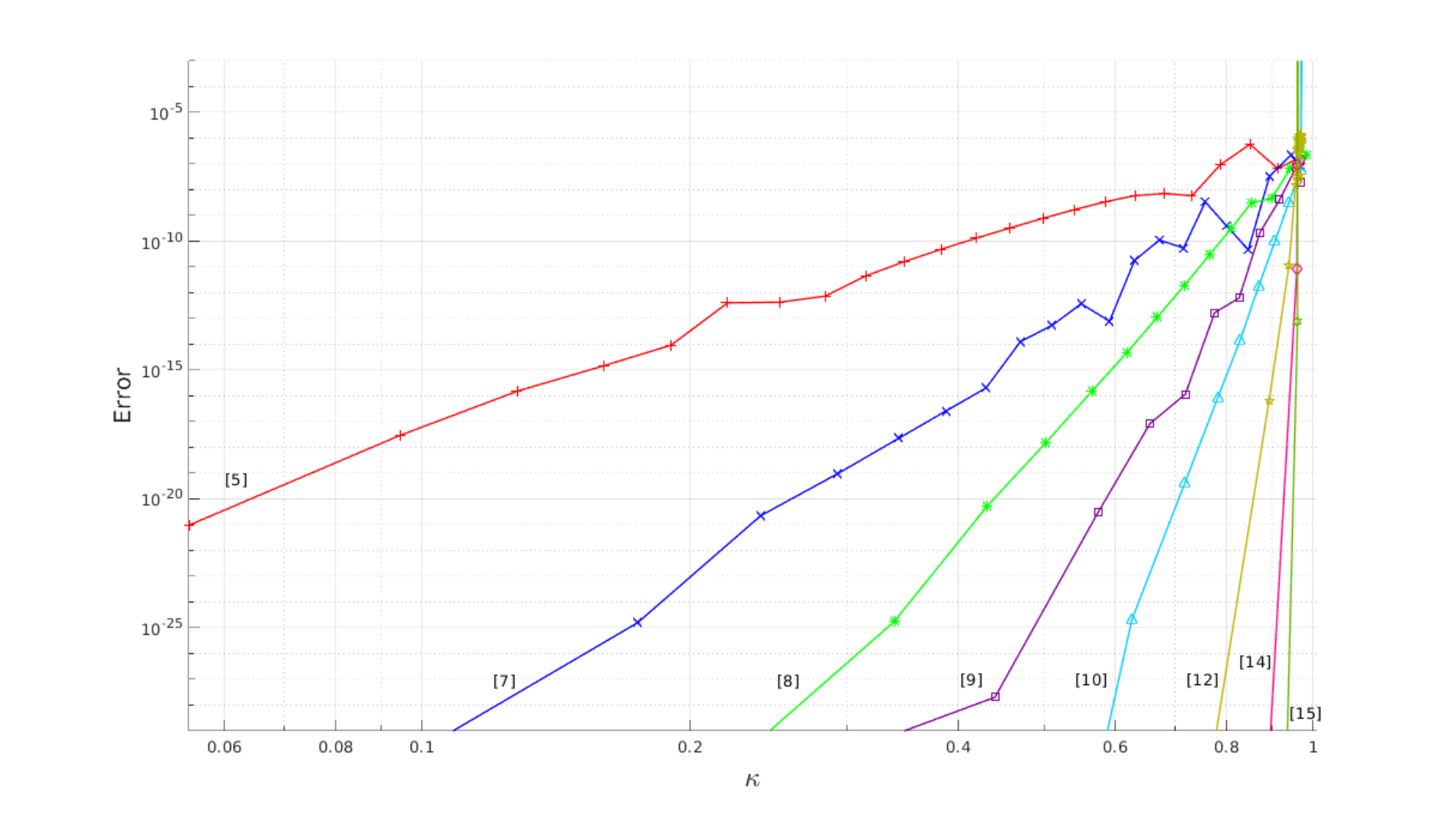} 
     \caption{Error evaluated as the minimum over $\theta$ of $\tilde{E}(\theta)$ in Eq.~(\ref{error_sob_1}) for periodic orbits in the Chirikov-Taylor map (\ref{stnmap_1}) as function of the parameter $\kappa$ using only the modified parameterization method (i.e., with no Newton-Gauss refinement) for different rotation numbers corresponding to Fibonacci ratios $[N]=F_{N-1}/F_{N}$ from $N=5$ to $N=15$.} 
     \label{error_fig_stn}
 \end{figure}

Figure \ref{error_fig_stn} (note the logarithmic scale) shows that, for a fixed period,  the error of the modified parameterization error exhibits a power law increase with $\kappa$, with exponent proportional to the period.  
This might be related to the observation that the modified parameterization method is attempting to approximate with a smooth continuous curve a periodic orbit that tends to have a fractal (self-similar) structure for large periods as 
$\kappa\rightarrow\kappa_c$. 
Consistent with this, for a fixed $\kappa$,
the error decreases as the period of the orbit increases.
It can also be appreciated from Fig.~\ref{error_fig_stn} that working in quadruple precision, the modified parameterization method gives  reasonable good estimates of the periodic orbits for values of $\kappa$ far from $\kappa_c$. The Newton-Gauss method may only be required, depending on the precision needs, for calculations
where the modulus of the residue of the periodic orbit is greater than $10^{-12}$.

\begin{figure}[h!]
     \centering
     \includegraphics[width=0.99\linewidth]{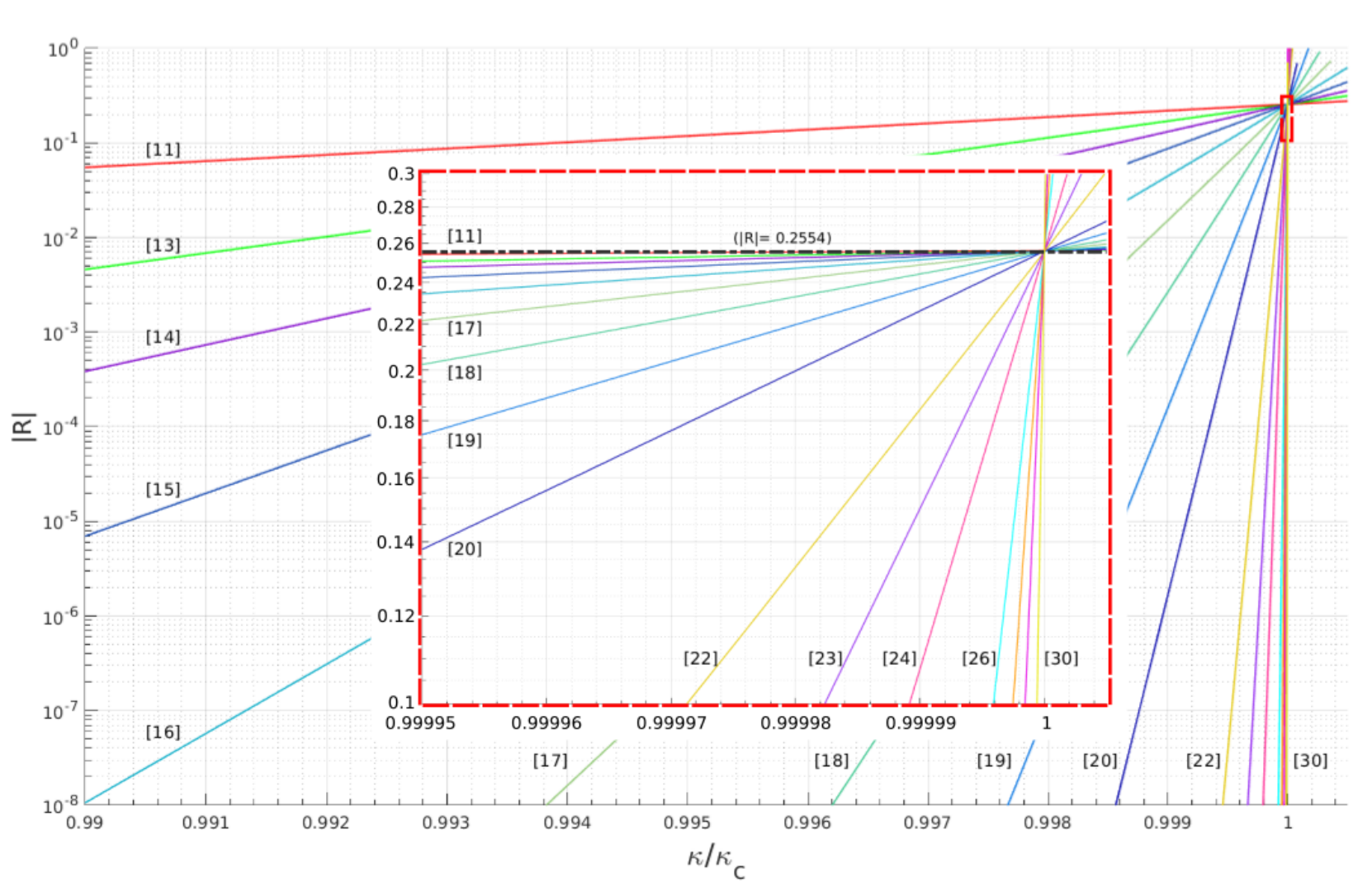}
     \caption{Absolute value of the residue of hyperbolic periodic orbits for the standard map (\ref{stnmap_1}) as function of the parameter $\kappa$ 
      computed with the compound method,
for rotation numbers corresponding to Fibonacci ratios $[N] \equiv F_{N-1}/F_{N}$ 
   from N=11 ($89/144$) to N=30 ($832040/1346269$). The zoom shown in the inset plot verifies the expected convergence of the residues to $|R|=0.2554$ at criticality $\kappa=\kappa_c$.    }
     \label{res_fig_stn}
 \end{figure}

The results of performing a numerical continuation using the compound method (i.e., the modified
parameterization method followed by the Newton-Gauss
refinement)  are shown in Table \ref{table_stn} and Fig.~\ref{res_fig_stn}. 
As expected, Fig.~\ref{res_fig_stn} shows that for a fixed $\kappa < \kappa_c$, the residue decreases when the period increases, which explains the decrease of the error with the increase of the period shown in  Fig.~\ref{error_fig_stn} and Table \ref{table_stn} because the  
accuracy of the parameterization method  is higher for small values of the residue. 
The inset in Fig.~\ref{res_fig_stn} shows that as expected, and consistent with well-known previous results, the residues converge to $|R|=0.2554$ at criticality $\kappa=\kappa_c$. 

%

\begin{table}[h!]
\renewcommand{\arraystretch}{1.6}
\centering
\begin{tabular}{c c c c} 
 \specialrule{.12em}{0em}{0em}
  $\Resize{0.45cm}{\frac{p}{q}}$   
      &  \begin{tabular}{c} Orbit \\ $x$ \end{tabular}
      &   ${E}$  & ${R}$ 
        \\ 
        \specialrule{.12em}{0em}{0em}
 {$\Resize{0.34cm}{\frac{5}{8}}$}  
   & $\Resize{5.8cm}{\phantom{oi}1.04262381263415710544\times 10^{-01}}$    
   & $\Resize{2.35cm}{1.7776\times 10^{-33}}$                               
   & $\Resize{2.6cm}{-2.3618\times 10^{-01}}$ \\                            
                 \specialrule{.04em}{0em}{0em}     
 $\Resize{0.52cm}{\frac{8}{13}}$   
    & $\Resize{5.8cm}{\phantom{oi}7.2292421040641022513\times 10^{-02}}$    
   & $\Resize{2.35cm}{3.8240\times 10^{-45}}$                               
   & $\Resize{2.6cm}{-2.1396\times 10^{-01}}$ \\                            
                 \specialrule{.04em}{0em}{0em}  
 $\Resize{0.52cm}{\frac{13}{21}}$ 
   & $\Resize{5.8cm}{\phantom{oi}5.06992260762055879603\times 10^{-02}}$    
   & $\Resize{2.35cm}{1.4227\times 10^{-30}}$                               
   & $\Resize{2.6cm}{-1.9779\times 10^{-01}}$ \\                            
                 \specialrule{.04em}{0em}{0em}  
 $\Resize{0.52cm}{\frac{21}{34}}$ 
    & $\Resize{5.8cm}{\phantom{oi}3.5022447826195364453\times 10^{-02}}$    
   & $\Resize{2.35cm}{3.0469\times 10^{-30}}$                               
   & $\Resize{2.6cm}{-1.6521\times 10^{-01}}$ \\                            
                 \specialrule{.04em}{0em}{0em}   
 $\Resize{0.52cm}{\frac{34}{55}}$  
    & $\Resize{5.8cm}{\phantom{oi}2.4048789198498598295\times 10^{-02}}$    
   & $\Resize{2.35cm}{1.5432\times 10^{-29}}$                               
   & $\Resize{2.6cm}{-1.2748\times 10^{-01}}$ \\                            
                 \specialrule{.04em}{0em}{0em}  
 {$\Resize{0.52cm}{\frac{55}{89}}$} 
   & $\Resize{5.8cm}{\phantom{oi}1.62171825022147208409\times 10^{-02}}$    
   & $\Resize{2.35cm}{2.8721\times 10^{-32}}$                               
   & $\Resize{2.6cm}{-8.2110\times 10^{-02}}$ \\                            
                 \specialrule{.04em}{0em}{0em}  
 $\Resize{0.67cm}{\frac{89}{144}}$ 
    & $\Resize{5.8cm}{\phantom{oi}1.0731773004662365072\times 10^{-02}}$    
   & $\Resize{2.35cm}{6.7440\times 10^{-52}}$                               
   & $\Resize{2.6cm}{-4.0823\times 10^{-02}}$ \\                            
                 \specialrule{.04em}{0em}{0em}  
 $\Resize{0.67cm}{\frac{144}{233}}$
    & $\Resize{5.8cm}{\phantom{oi}6.9344882456917456713\times 10^{-03}}$    
   & $\Resize{2.35cm}{5.2130\times 10^{-55}}$                               
   & $\Resize{2.6cm}{-1.3101\times 10^{-02}}$ \\                            
                 \specialrule{.04em}{0em}{0em}  
 $\Resize{0.67cm}{\frac{233}{377}}$ 
   & $\Resize{5.8cm}{\phantom{oi}4.39068507873134101191\times 10^{-03}}$    
   & $\Resize{2.35cm}{1.9267\times 10^{-32}}$                               
   & $\Resize{2.6cm}{-2.0969\times 10^{-03}}$ \\                            
                 \specialrule{.04em}{0em}{0em} 
 $\Resize{0.67cm}{\frac{377}{610}}$ 
   & $\Resize{5.8cm}{\phantom{oi}2.74217144026659842688\times 10^{-03}}$    
   & $\Resize{2.35cm}{9.5578\times 10^{-30}}$                               
   & $\Resize{2.6cm}{-1.0800\times 10^{-04}}$ \\                            
                 \specialrule{.04em}{0em}{0em}  
 $\Resize{0.67cm}{\frac{610}{987}}$ 
   & $\Resize{5.8cm}{\phantom{oi}1.70222914399736324091\times 10^{-03}}$    
   & $\Resize{2.35cm}{9.8163\times 10^{-45}}$                               
   & $\Resize{2.6cm}{-8.9147\times 10^{-07}}$ \\                            
                 \specialrule{.04em}{0em}{0em}
  $\Resize{0.73cm}{\frac{987}{1597}}$
   & $\Resize{5.8cm}{\phantom{oi}1.05402094177906953965\times 10^{-03}}$    
   & $\Resize{2.35cm}{6.5275\times 10^{-46}}$                               
   & $\Resize{2.6cm}{-3.7914\times 10^{-10}}$ \\                            
                  \specialrule{.12em}{0em}{0em}       
\end{tabular}
\renewcommand{\arraystretch}{1}
\caption{Computational error, $E$ in (\ref{E_NG}), residue, $R$ in (\ref{greene_res}), and $x$-component of hyperbolic periodic orbit location closest to 
$x=0$ for the first few $p/q$ inverse golden mean approximants  in the Chirikov-Taylor map (\ref{stnmap_1}) with $\kappa=0.9600$, computed using the 
compound method.
}
\label{table_stn}
\end{table}

	\subsection{Rational harmonic map (RHM)}
\begin{figure}[h!]
 \centering
   \includegraphics[width=0.95\linewidth]{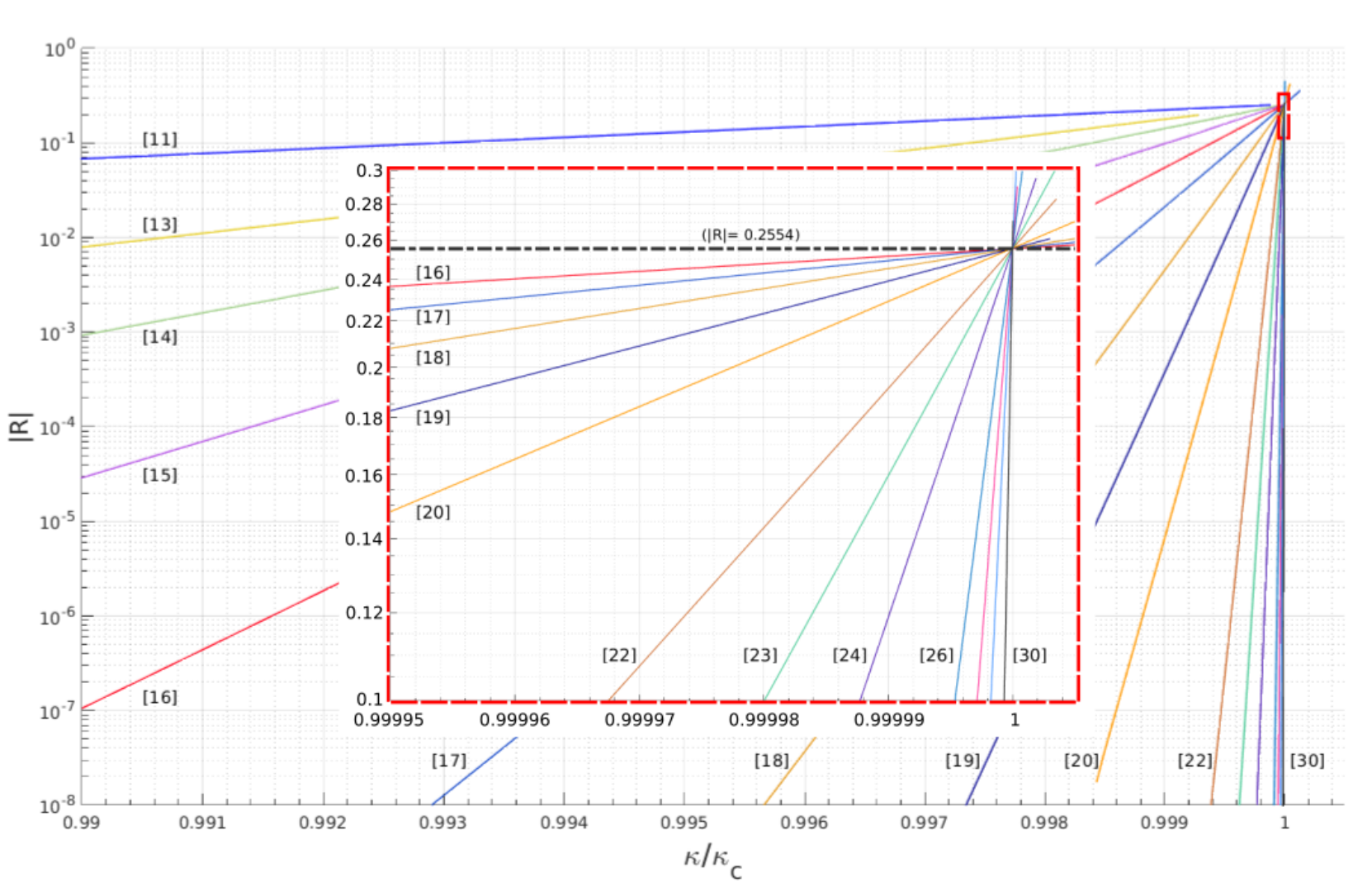} 
   \caption{Amplitude of ${R}$ of the approximated hyperbolic points for the RHM (\ref{stnmap_0}) and (\ref{anmap_func}) as function of the parameter $\kappa$ 
     and $(\alpha,\beta)=(3.0,0.4)$ computed with the compound method,
for different rotation numbers corresponding to Fibonacci ratios $[N] \equiv F_{N-1}/F_{N}$ 
   from N=11 ($89/144$) to N=30 ($832040/1346269$).}
     \label{res_fig_rhm}
\end{figure}

\begin{figure}[h!]
   \centering
   \includegraphics[width=0.99\linewidth]{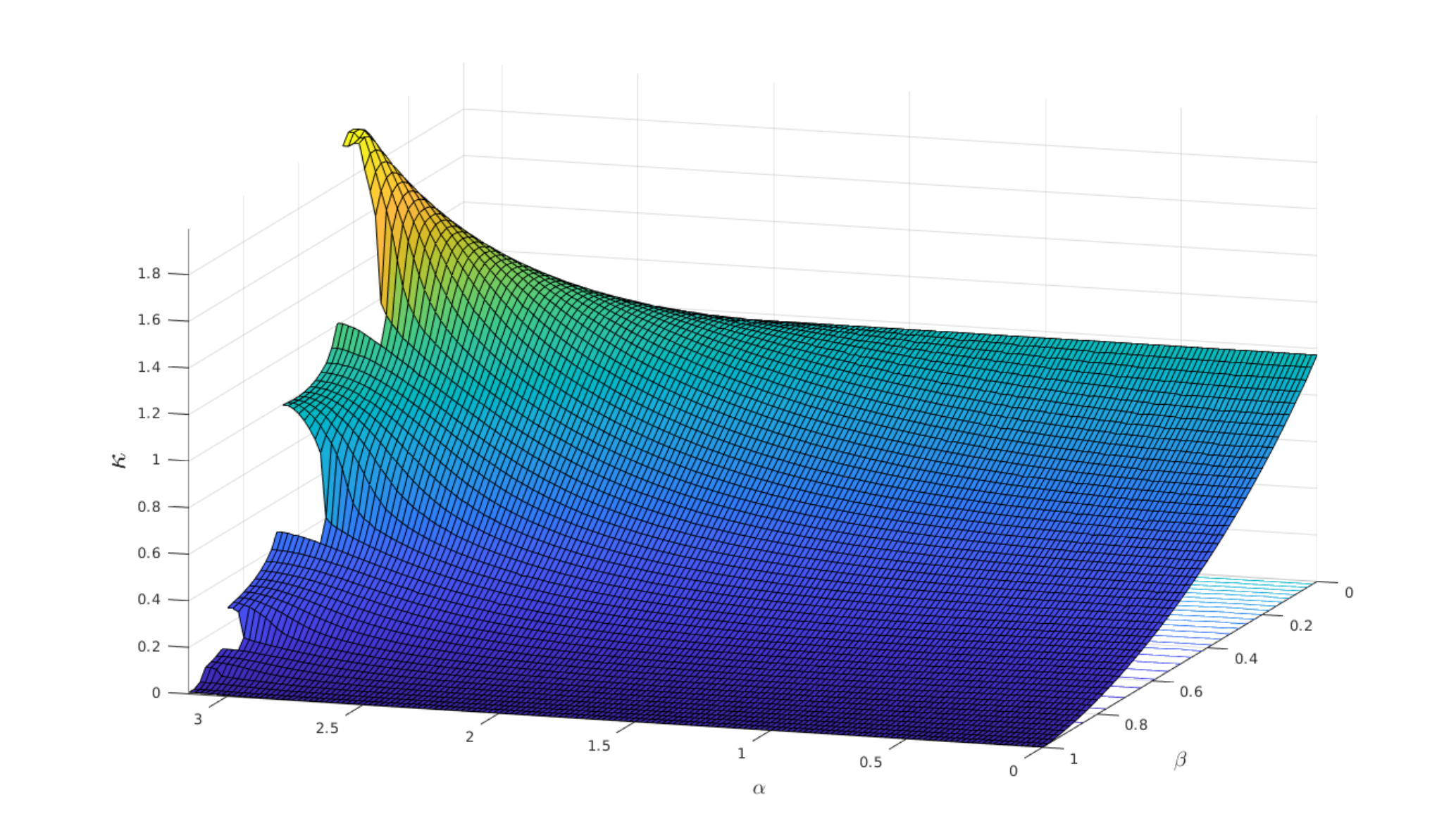}   
   \caption{\color{blue}Criticality isosurface of the rational harmonic map (\ref{anmap_func}) on the parameter space, calculated for the 
     $4181/6765$-periodic orbit ($N=19$) 
      for the critical residue value $|R| = 0.2545$ estimated in Fig.~\ref{res_fig_rhm}.}%
     \label{anmap_par}
 \end{figure}

\begin{figure}[h!]
   \centering  
   \includegraphics[width=0.99\linewidth, trim = {0.8cm 0.0cm 0.8cm 0.5cm}]{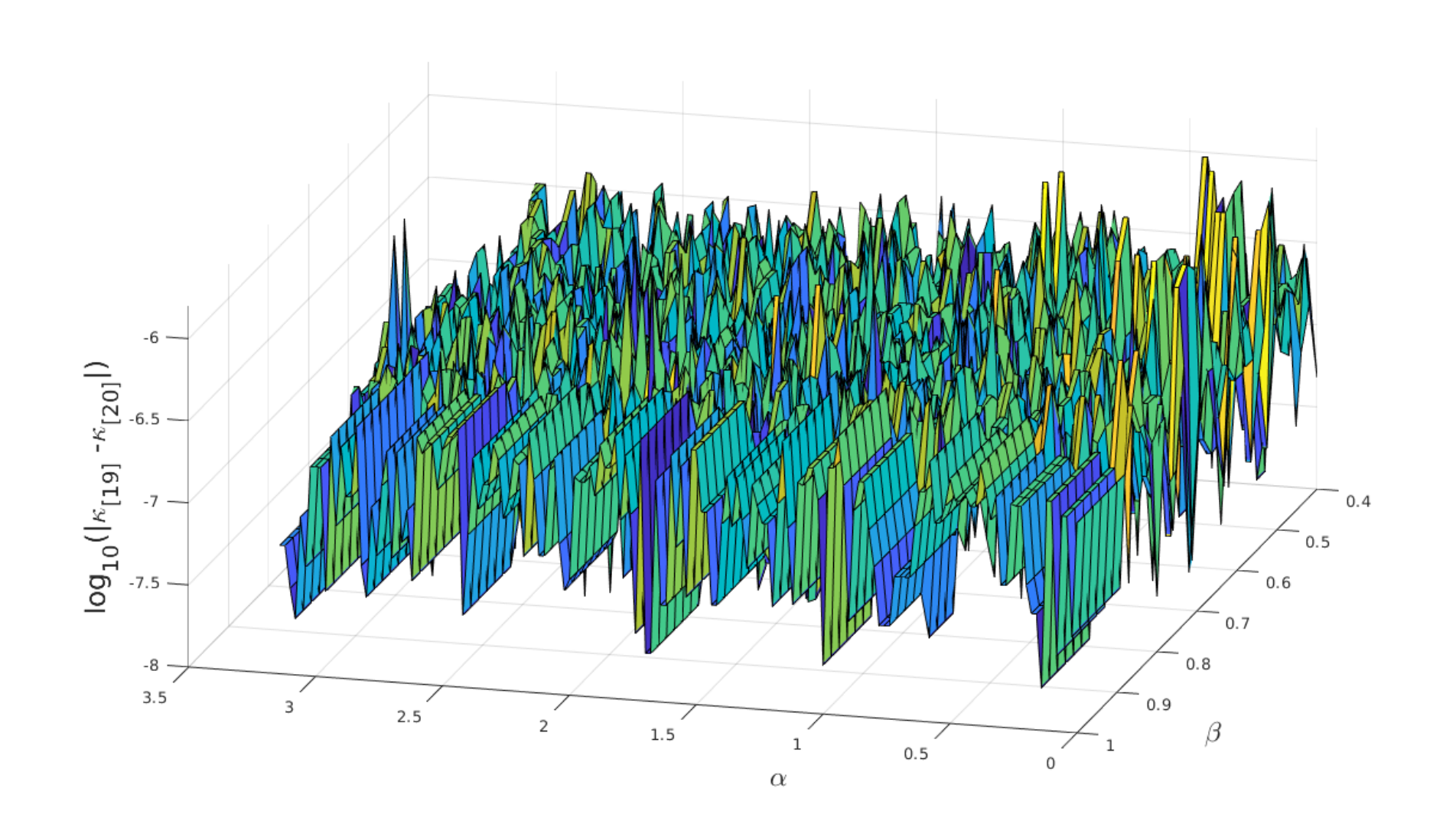}    
   \caption{
     Absolute difference between 
     two manifolds on the parameter space corresponding to: 
     $[19]= 4181/6765$- and $[20]=6765/10946$- periodic orbits of the rational harmonic map  
     (\ref{anmap_func}) obtained by continuation over parameter $\kappa$ for a regular grid on $(p_a,p_m)$, starting from the integrable case $\kappa=0$ and stopping for $|R| = 0.2554$.}
     \label{anmap_par_dif}
 \end{figure}

{\color{blue} As an example illustrating the application of the compound method to a map without a known useful involution decomposition, we consider the rational harmonic map defined in (\ref{stnmap_0}) and (\ref{anmap_func}).  This map has three independent parameters and, as a starting point, we fix 
$(\alpha,\beta)=(3.0,0.4)$ and study the behavior of the residues as function of $\kappa$, 
thus reducing the computation to a single-parameter scan. 
The numerical continuation for this map in Ref.~\cite{Petrov_Olvera08} was found to be very sensitive to the variation of parameters $\alpha$ and $\beta$.

The main results of this numerical study are summarized in Fig.~\ref{res_fig_rhm} that shows the absolute value of the residues of the hyperbolic orbits with periods, $F_{N-1}/F_N$, where $\{F_N\}$ are the Fibonacci convergents of the inverse golden mean. As expected, the residues  increase with $\kappa$. But, most importantly, at $\kappa=1.73360453\ldots$ all the residues are observed to converge  to 
$|R|=0.2554 \ldots$\footnote{\color{blue}  Which coincides up to 4 digits of the value reported in the literature, e.g. Ref.~\cite{Reichl04}}. 
Within numerical error, this value agrees with the one identified since the pioneering work Refs.~\cite{Greene79,Kadanoff81,shenker82,mackay1983renormalization,MacKay82} as the universal critical residue for reversible twist maps along the dominant symmetry line. 
Our findings provide numerical evidence of the universality of the critical residue beyond the realm of reversible maps. 
As previously done in Table \ref{table_stn} for the standard map, Table \ref{table_analy} presents, for the rational harmonic map, the values of the angular position of several periodic orbits to appreciate the degree of approximation in the Modified Parameterization Method and the Newton-Gauss method.

As a next step we study the criticality of the rational harmonic map in the 3-dimensional 
$(\alpha,\beta,\kappa)$ space. Note that using the symmetry $(x,y;\kappa,\alpha,\beta)\mapsto (x-1/2,y;\kappa,\alpha+\pi,-\beta)$ it is enough to consider the $\mathbb{R}^2 \times \mathbb{R}^+$ subspace.  
To guide this study, based on the results in Fig.~\ref{res_fig_rhm}, we make the ansatz that,
for generic parameter values, the critical residue remains $|R_c|=0.2554 \ldots$ and search for the 2-dimensional critical manifold 
defined as the locus satisfying  $R_{[N]}(\alpha,\beta,\kappa)=R_c$, where $R_{[N]}$ is the absolute value of the residue of the $F_{N-1}/F_N$ hyperbolic periodic orbit for parameters $(\alpha,\beta,\kappa)$.

We computed critical iso-surfaces corresponding to $F_{10}/F_{11}=89/144$, $F_{17}/F_{18}=2584/4181$, $F_{18}/F_{19}=4181/6765$ and $F_{19}/F_{20}=6765/10946$. As an example, 
Figure \ref{anmap_par} shows the case $R_{[19]}(\alpha,\beta,\kappa)=R_c$
corresponding to $F_{18}/F_{19}=4181/6765$.
The iso-surface exhibits smooth behavior, except near the $\alpha=\pi$ boundary where  it shows rapid variations and ``cusps" in $\kappa$ as $\beta$ changes which suggest the possibility that the surface might have folds  in the region $\alpha \sim \pi$.  
To numerically construct this critical  iso-surface we setup a uniform $[0,\pi) \times[0,1)$ grid in the  $(\alpha,\beta)$ plane and performed a continuation over the parameter $\kappa$ for each grid point. We  use the modified parameterization method (MPM) to first continue the orbits from the integrable case $\kappa=0$ until a certain $\kappa_1$ value for which the residues of the orbits could be large enough to use the Newton-Gauss method to reliable and more accurately continue the orbits to a critical value of parameter $\kappa$ and beyond.
The stopping criteria for the continuation is that the residue of the orbit outgrows the critical threshold$|R| = 0.2554$
or that the error $E_1$ can not be reduced by the Newton step of the parameterization method below the a priori upper bound $E_1\lesssim 10^{\mathbf{-18}}$, 
although the Newton-Gauss method reduces it to less than $10^{-30}$.

To test the validity of the renormalization  ansatz (universality of critical residue) we computed the critical manifold for the next convergent,
$F_{19}/F_{20}=6765/10946$, and found that the iso-surfaces 
$R_{[11]}(\alpha,\beta,\kappa)=R_c$, $R_{[18]}(\alpha,\beta,\kappa)=R_c$,
$R_{[19]}(\alpha,\beta,\kappa)=R_c$ and $R_{[20]}(\alpha,\beta,\kappa)=R_c$ to be practically indistinguishable. 
For the range of  parameters considered (which exclude the $\alpha \sim \pi$ region)  these iso-surfaces can be represented as the graph of functions of the form $\kappa_{[N]}=\kappa_{[N]}(\alpha,\beta)$ where $R_{[N]}(\alpha,\beta,\kappa_{[N]})=R_c$. Based on this, we quantify the proximity between these surfaces using the difference in the height of the graphs. As Fig.~\ref{anmap_par_dif} 
shows (note the logarithmic scale), the difference in the heights, $\kappa_{[20]}(\alpha,\beta)-\kappa_{[19]}(\alpha,\beta)$, is in the range $(10^{-8},10^{-6})$, 
and a similar behavior was found for the differences between the other surfaces, i.e., $\kappa_{[N]}(\alpha,\beta)-\kappa_{[M]}(\alpha,\beta)$.
These results provide numerical evidence of the validity of the  ansatz and lend further support to the universality of the $R_c=0.2554 \ldots$ convergence at criticality which is a cornerstone of the renormalization theory of the destruction of invariant circles in twist maps.

Beyond serving as an illustration of the accuracy and efficiency of the proposed compound method, these results provide, in our understanding for the first time, numerical evidence of the universality 
(in terms of the residue convergence) of the breakup of the inverse golden mean invariant circle for 
multiparametric
maps lacking involution decomposition and symmetry lines.}

It can be seen from Fig.~\ref{anmap_par} that the behavior of the critical value of $\kappa$ as function of $\beta$ seems to agree with what expected from the perturbation function (\ref{anmap_func}), bigger values of $\beta$ should yield lower values of $\kappa_c$. 
It is however unexpected to find that for $\alpha \sim \pi$ there appear irregularities that suggest that the critical surface may have folds and also, in the same region, $\kappa$ rise to values bigger than the critical $\kappa_c$, the analytic critical value for $\beta=0$ and any $\alpha$. Because of this peculiar behavior of the manifold, it was chosen to compute periodic orbits for $(\alpha,\beta)=(3.0,0.4)$, close to a fold but yet not the reversible case of $\alpha=\pi$. Some of the computed orbits are presented in table \ref{table_analy}.
%
\begin{table}[h!]
\centering
\begin{tabular}{c g g g} 
 \specialrule{.12em}{0em}{0em}
\rowcolor{white}
  $\Resize{0.4cm}{\frac{p}{q}}$   
      &  \begin{tabular}{c} Orbit \\ $x$ \end{tabular}
      &   ${E}$ & ${R}$
      \\     \specialrule{.12em}{0em}{0em}  
 \rowcolor{white}
 \multirow{2}{*}{$\Resize{0.52cm}{\frac{8}{13}}$}    
   & $\Resize{5.6cm}{9.84400324333751592206\times 10^{-01}}$           
   & $\Resize{2.51cm}{1.1780\times 10^{-05}}$             
    &$\Resize{2.8cm}{-2.1690\times 10^{-01}}$ \\                       
   & $\Resize{5.6cm}{9.84194430418146115882\times 10^{-01}}$           
   & $\Resize{2.51cm}{4.8935\times 10^{-39}}$                          
    &$\Resize{2.8cm}{-2.1777\times 10^{-01}}$ \\                       
   \specialrule{.04em}{0em}{0em}   
 \rowcolor{white}
 \multirow{2}{*}{$\Resize{0.52cm}{\frac{21}{34}}$} 
   & $\Resize{5.6cm}{9.97753743054888484767\times 10^{-01}}$           
   &$\Resize{2.51cm}{4.0668\times 10^{-07}}$            
    &$\Resize{2.8cm}{-1.8059\times 10^{-01}}$ \\                       
   & $\Resize{5.6cm}{9.97722363363172234098\times 10^{-01}}$           
   & $\Resize{2.51cm}{1.6284\times 10^{-47}}$                          
    &$\Resize{2.8cm}{-1.8095\times 10^{-01}}$ \\                       
   \specialrule{.04em}{0em}{0em}     
 \rowcolor{white}
 \multirow{2}{*}{$\Resize{0.52cm}{\frac{55}{89}}$} 
   & $\Resize{5.6cm}{4.35066410514029663774\times 10^{-04}}$           
   & $\Resize{2.51cm}{3.1674\times 10^{-05}}$            
    &$\Resize{2.8cm}{-2.8817\times 10^{-01}}$ \\                       
   & $\Resize{5.6cm}{4.42076348408883988105\times 10^{-04}}$           
   & $\Resize{2.51cm}{1.8242\times 10^{-37}}$                          
    &$\Resize{2.8cm}{-1.0589\times 10^{-01}}$ \\                       
   \specialrule{.04em}{0em}{0em} 
 \rowcolor{white}
 \multirow{2}{*}{$\Resize{0.67cm}{\frac{144}{223}}$} 
   & $\Resize{5.6cm}{4.09885888135278521212\times 10^{-04}}$           
   & $\Resize{2.51cm}{3.5713\times 10^{-13}}$             
    &$\Resize{2.8cm}{-1.3204\times 10^{-01}}$ \\                       
   & $\Resize{5.6cm}{4.54519288406945185848\times 10^{-04}}$           
   & $\Resize{2.51cm}{8.2324\times 10^{-29}}$                          
    &$\Resize{2.8cm}{-2.5599\times 10^{-02}}$ \\                       
   \specialrule{.04em}{0em}{0em} 
 \rowcolor{white}
 \multirow{2}{*}{$\Resize{0.67cm}{\frac{377}{610}}$} 
   & $\Resize{5.6cm}{4.54770852734730605323\times 10^{-04}}$           
   & $\Resize{2.51cm}{1.4142\times 10^{-07}}$             
    &$\Resize{2.8cm}{-2.5479\times 10^{-02}}$ \\                       
   & $\Resize{5.6cm}{5.10890552501945729879\times 10^{-04}}$           
   & $\Resize{2.51cm}{1.2373\times 10^{-29}}$                          
    &$\Resize{2.8cm}{-6.2363\times 10^{-04}}$ \\                       
    \specialrule{.04em}{0em}{0em} 
 \rowcolor{white}
 \multirow{2}{*}{$\Resize{0.76cm}{\frac{987}{1597}}$} 
   & $\Resize{5.6cm}{2.18102889348964055253\times 10^{-04}}$           
   & $\Resize{2.51cm}{5.9188\times 10^{-14}}$             
    &$\Resize{2.8cm}{-9.2787\times 10^{-08}}$ \\                       
   & $\Resize{5.6cm}{2.88969724223810721485\times 10^{-04}}$           
   & $\Resize{2.51cm}{9.5395\times 10^{-31}}$                          
    &$\Resize{2.8cm}{-4.0099\times 10^{-08}}$ \\                       
    \specialrule{.04em}{0em}{0em}      
 \rowcolor{white}
 \multirow{2}{*}{$\Resize{0.76cm}{\frac{2584}{4181}}$} 
   & $\Resize{5.6cm}{5.70931653171076079438\times 10^{-03}}$           
   &$\Resize{2.51cm}{4.4036\times 10^{-28}}$             
    &$\Resize{2.8cm}{-3.2521\times 10^{-19}}$ \\                       
   & $\Resize{5.6cm}{5.70811191353601822914\times 10^{-03}}$           
   & $\Resize{2.51cm}{1.1872\times 10^{-31}}$                          
    &$\Resize{2.8cm}{-3.2944\times 10^{-19}}$ \\                       
     \specialrule{.04em}{0em}{0em}    
 \rowcolor{white}
 \multirow{2}{*}{$\Resize{0.82cm}{\frac{6765}{10946}}$} 
   & $\Resize{5.6cm}{1.236939703348994242120533777\times 10^{-02}}$    
   &$\Resize{2.51cm}{2.8990\times 10^{-32}}$            
    &$\Resize{2.8cm}{-8.6237\times 10^{-24}}$ \\                       
   & $\Resize{5.6cm}{1.236939703348994242120533769\times 10^{-02}}$    
   & $\Resize{2.51cm}{9.5204\times 10^{-29}}$                          
    &$\Resize{2.8cm}{-2.0985\times 10^{-24}}$ \\                       
   \specialrule{.12em}{0em}{0em} 
\end{tabular}
\caption{
Computational error, $E$, residue, $R$, and $x$-component of hyperbolic periodic orbit location closest to 
$x=0$ for the first few $p/q$ inverse golden mean approximants  in the rational harmonic map (\ref{anmap_func}) with $(\kappa,\alpha,\beta)=(1.7150,3.0,0.4)$, computed using the modified parameterization method (white rows) and the compound method (gray rows). 
For the modified parameterization method $E$ and $R$ are defined in  (\ref{error_sob_1}) and (\ref{aprox_res}), and for the compound method the definitions are (\ref{E_NG}) and (\ref{greene_res}).}
\label{table_analy}
\end{table}

\subsection{Beyond the break up of invariant circles, reoccurrence, and other results}
\begin{figure}[h!]
     \centering  
     \includegraphics[width=0.98\linewidth]{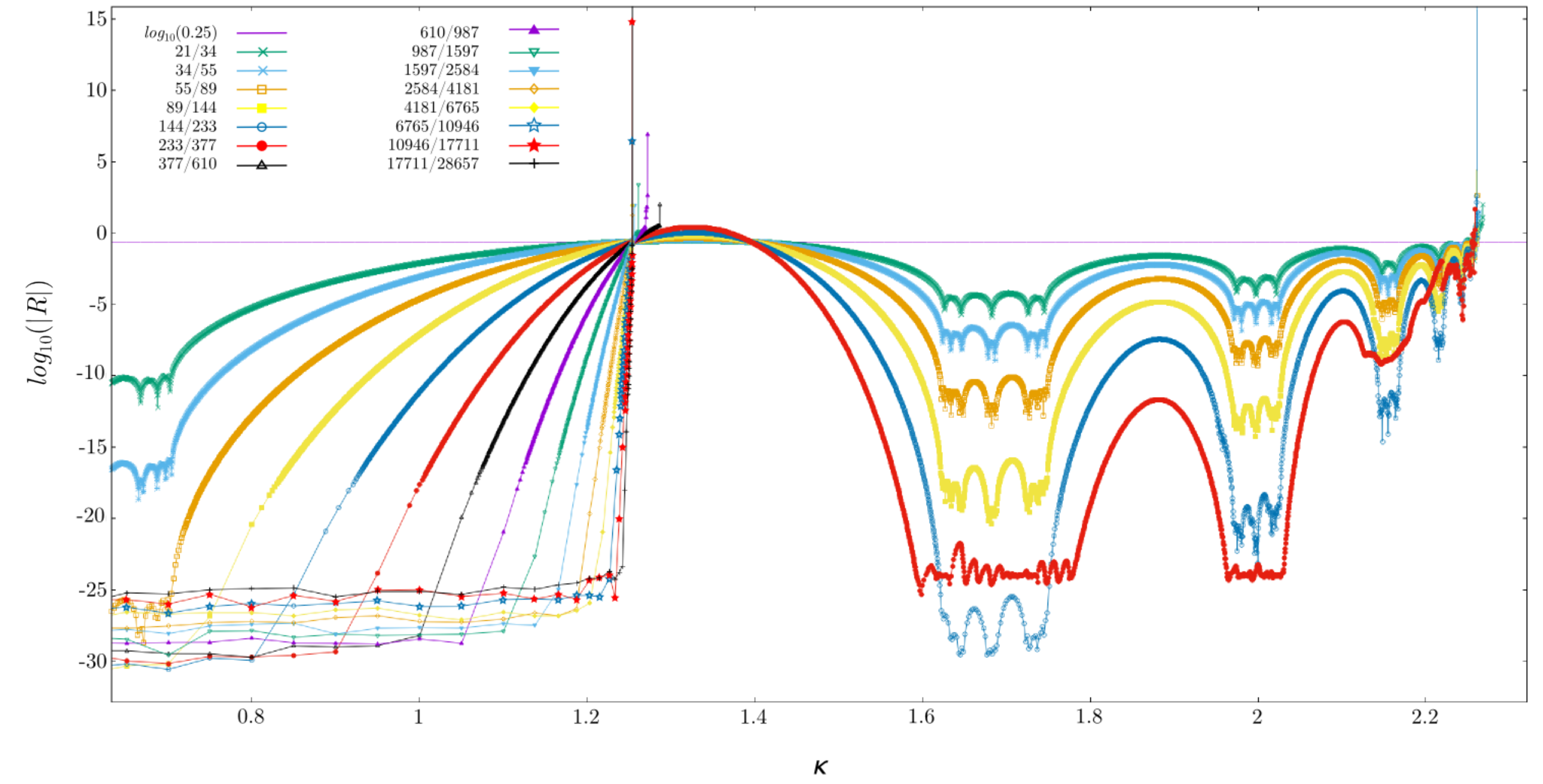}  
     \caption{Residue of periodic orbits as function of the parameter $\kappa$ in the rational harmonic map 
     (\ref{anmap_func}) for $(\alpha,\beta)=(3.0,0.4)$.}
     \label{perio_resurgence}
 \end{figure}

A surprising result, not originally envisioned in the implementation of the proposed compound  method, is that it is possible to continue the computation of periodic orbits beyond the break up of the corresponding invariant circle. It is increasingly difficult for high order periodic orbits, but it is possible to compute them and obtain particular information about the {\color{blue}\emph{reoccurrence}} of a given invariant circle. This {\color{blue}\emph{reoccurrence}} can be understood if we consider that the continuation on parameter $\kappa$ is done along vertical lines that eventually may cross through folds in the critical manifold. 
The {\color{blue}\emph{reoccurrence}} can also be appreciated in the residue versus $\kappa$ plots in Fig.~\ref{perio_resurgence} 
where for $(\alpha,\beta)=(3.0,0.4)$ in the rational harmonic map the residues of five different $F_{n-1}/F_n$-periodic orbits beyond the critical value for the break up,  drop below $|R|\lesssim 0.25$ for $\kappa\in[1.41,2.25]$, which signals the {\color{blue}\emph{reoccurrence}} of the invariant circle. 

\begin{figure}[h!]
 \centering
   \includegraphics[width=0.99\linewidth]{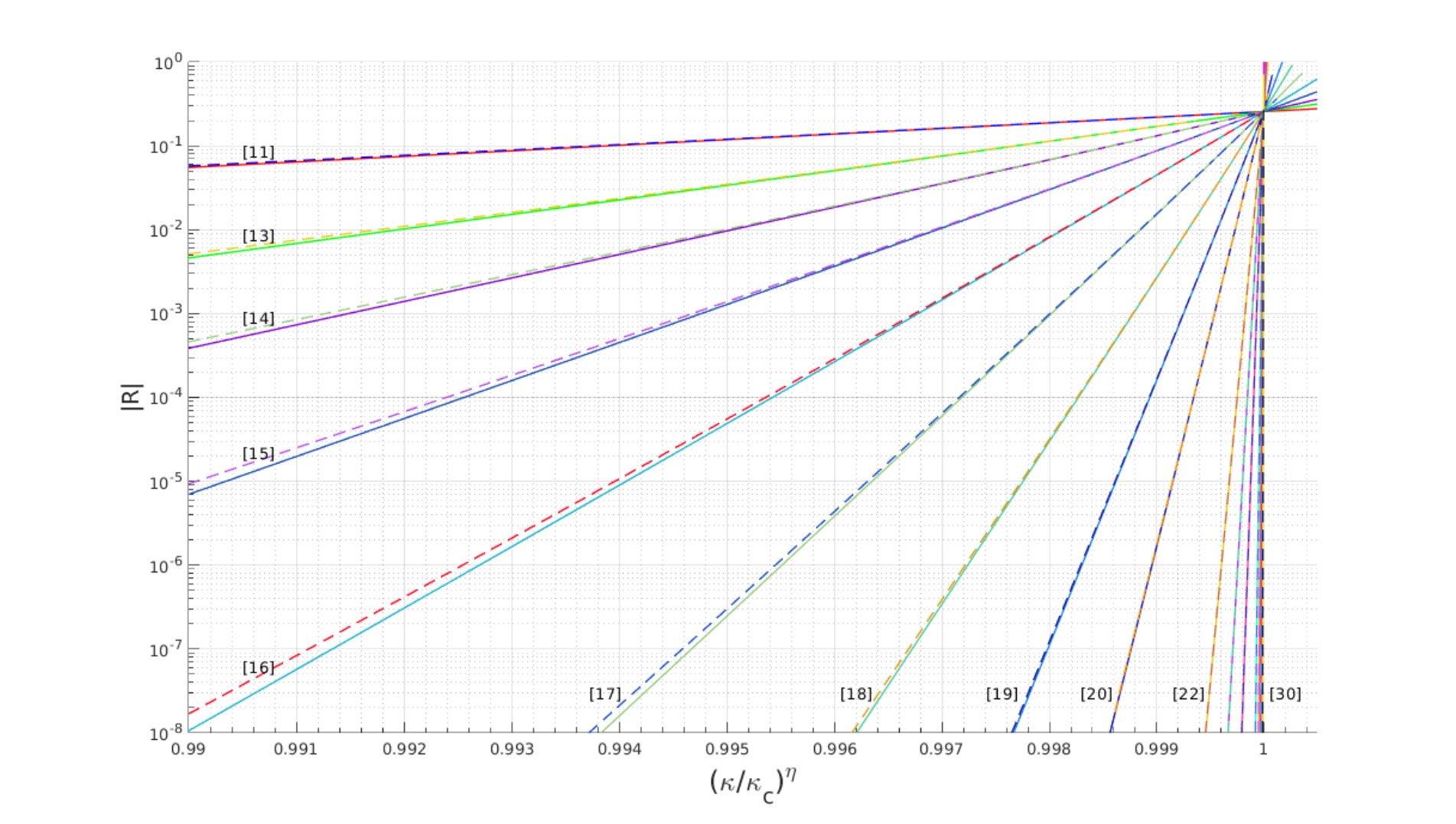} 
   \caption{Comparison between the residues $R$ of hyperbolic periodic orbits for the standard map (continuous lines) and the rational harmonic map (dashed lines) (\ref{anmap_func}) for $(\alpha,\beta)=(3.0,0.4)$, both as function of a scaled parameter $(\kappa/\kappa_c)^\eta$. Where $(\kappa_c,\eta)=(0.9716354,1.0)$ for the standard map and $(\kappa_c,\eta)=(1.73360453,0.885)$ for the RHM.
   The orbits rotation numbers correspond to the Fibonacci ratios $[N] \equiv F_{N-1}/F_{N}$ from N=11 ($89/144$) to N=30 ($832040/1346269$). }
     \label{res_comp_NG}
\end{figure}

{\color{blue} Another interesting feature, found when comparing the results from standard map in Fig.~\ref{res_fig_stn} and the rational harmonic map in Fig.~\ref{res_fig_rhm} for $R_{[N]}=(3.0,0.4,\kappa)$, is that, as shown in Fig.~\ref{res_comp_NG}, rescaling the horizontal axis as $(\kappa/\kappa_c)^\eta$, where $\eta$ is numerically determined exponent, reveals a correlation between the values of $R_{[N]}$ in both maps near the critical value $\kappa_c$.}

\section{\large{Discussion}}
\label{sec_Discussion}
The method proposed in this work can be applied to a wide variety of maps for the efficient and accurate computation of  periodic orbits without the use of symmetries to simplify the search. 
In the cases considered, the only requirements were the symplectic and twist conditions needed to justify some hypothesis of the parameterization method. The method as presented does not include the non-twist case when the periodic orbits are close to a shearless invariant circle, see Refs.~\cite{del1994dynamics,fuchss06}. 
However, in principle the method can be generalized to any map for which the standard parameterization method has been, or can be, applied (see Ref.~\cite{Haro16}), so an implementation for non-twist maps is not only possible but work in progress.

It was found that the modified parameterization method is efficient and very fast to find periodic orbits in cases where the parameters values are relatively close to the integrable case ($R \lesssim 10^{-8}$). 
On the other hand, the Newton-Gauss method was found to be very efficient for finding periodic orbits  close to the critical parameter values ($R\gtrsim 10^{-8}$), although it is relatively slower as continuation method for parameter values farther of criticality.
Both methods can be used to find periodic orbits but each one excels in complementary regimen.  The compound method takes the strengths of each scheme and allows the computation of periodic orbits of high period in a relatively fast and efficient way.

The application of these methods to non-autonomous systems is possible for the case of time-periodic maps as those discussed in Ref.~\cite{nasm_paper}. 
The composition of $q$ iterations of a given periodic non-autonomous twist map $T$ of period $q$, generally yields a very convoluted map $T^q$ that may not have obvious symmetries and in general there are regions in phase space where the map is non-twist. The compound method does not have such requirements so
it allows the computation of periodic orbits in a $T^q$ map that later can be related to the original map.
As described in Ref.~\cite{nasm_paper}, results from renormalization theory and
the computation of critical exponents only makes sense in the compound
map. In the non-autonomous maps, the phase space is not two dimensional and the orbit is not well defined in the projection to the cylinder.
The authors will present results for non autonomous periodic maps in a future publication.

The numerical experiments performed allowed us to verify results from the renormalization theory in two different ways. First, following the continuation of periodic orbits of rotation numbers $F_{N-1}/F_N$  approximating the inverse golden mean invariant circle an agreement between the critical residue for both the standard map and the rational harmonic map was found.
Secondly, to further explore the renormalization properties of the rational harmonic map for a wide range of parameters values we constructed the \textcolor{blue}{critical manifold} by performing a parameter scan in the $(\alpha, \beta)$ plane and computing, for each $\alpha$ and $\beta$, the value of $\kappa$ for which the absolute value of the  residue, $|R|$, reached the universal critical residue for twist maps.
\textcolor{blue}{The critical manifold was constructed for}
different periodic orbits with rotation numbers given by the Fibonacci's ratios approaching the inverse golden mean and 
\textcolor{blue}{a convergence of the successive critical manifolds was observed 
up to an error of order $10^{-5}$. That is,} 
 the critical parameter value $\kappa$ \textcolor{blue}{found} for every $(\alpha,\beta)$ pair was almost the same, as predicted by renormalization theory.
Finally, the {\color{blue}\emph{reoccurrence}} of the inverse golden mean invariant circle for an extended value of the parameter $\kappa$ was observed.  

To the authors knowledge, a verification of renormalization theory 
prediction of the residues' convergence at criticality for maps
for maps without
symmetries, like the RHM, has not been performed in the past.
The results show the consistency between Greene's residue conjecture and the universality of
critical scalings predicted by renormalization theory.
We believe that the use of Greene's residue criterion together with the continuation method
we propose here, may help to determine when a broken invariant circle may reappear for
larger values of the parameters after it has broken up.

The extension of the proposed modified parameterization method to higher dimensions is only limited by the intrinsic restrictions of the parameterization method, which is easier to extend  in even dimensions. The method might also be extended to  systems with small dissipation  \cite{Cal-Cel-Lla-13}.
On the other hand,  the extension of the Newton-Gauss method to higher dimensions seems plausible, with 
the only limitation that the size of the matrix $DG(\mathbf{z})$ and the number of operations will naturally increase. We should emphasize that in both cases the increases will remain related to the order of the periodic orbit in an analogous manner to the 2-D case.
As a particular example, we are currently working on the implementation of the method to high-dimensional Froeschl\'e type maps.

\section*{\large{Acknowledgments}}
This work was founded by PAPIIT IN112920, IN110317, IA102818 and IN101020, FENOMEC-UNAM and by the Office of Fusion Energy Sciences of the US Department
of Energy at Oak Ridge National Laboratory, managed by UT-Battelle, LLC, for the U.S.Department of Energy under contract DE-AC05-00OR22725. 
This material is based upon work supported by the National Science Foundation under Grant No. DMS-1440140 while R.C. and D.M. were in residence at the Mathematical Sciences Research Institute in Berkeley, California, during the Fall 2018 semester.
We are indebted with \`A. Haro for sharing with us the ideas behind the implementation of the Newton-Gauss method. Also it is a pleasure to acknowledge the insightful discussions we had at Berkeley with J.D. Meiss, N. Petrov and A. Wurm, as well as the valuable feedback from the referees. %
We also express our gratitude  to the graduate program in Mathematics of UNAM for making the GPU servers available to perform our computations and especially to Ana Perez for her invaluable help.

\bibliographystyle{unsrt}

\end{document}